
 \input amstex
 \documentstyle{amsppt}
 \nologo \NoBlackBoxes \magnification=\magstep1
 \pageheight{45pc}
 \topmatter
 \title  On Cohomology of the Square of an Ideal Sheaf \endtitle
 \affil University of North Carolina \endaffil
 \author Jonathan Wahl \endauthor
 \address Department of Mathematics,
   University of North Carolina,
   Chapel Hill, NC 27599-3250 \endaddress
 \email jw\@math.unc.edu \endemail
 \abstract
For a smooth subvariety $X\subset\Bbb P^N$, consider (analogously to
projective normality) the vanishing condition $H^1(\Bbb P^N,\Cal
I^2_X(k))=0$, $k\ge3$.  This condition is shown to be satisfied for all
sufficiently large embeddings of a given $X$, and for a
Veronese embedding of $\Bbb P^n$.

For $C\subset\Bbb P^{g-1}$, the canonical embedding of a non-hyperelliptic
curve, this condition guarantees the vanishing of some obstruction groups
to deformations of the cone. Recall that the tangents to deformations are
dual to the cokernel of the Gaussian-Wahl map.

\proclaim{Theorem} Suppose the Gaussian-Wahl map of $C$ is not surjective
and the vanishing condition is fulfilled. Then $C$ is {\bf extendable}:
 it is a
hyperplane section of a surface in $\Bbb P^g$ not the cone over $C$.
\endproclaim

Such a surface is a K3 if smooth, but it could have serious
singularities.

\proclaim{Theorem} For a general curve of
genus $\ge3$, this vanishing holds. \endproclaim

\proclaim{Conjecture} If the Clifford index is $\ge3$, this vanishing holds.
\endproclaim
\endabstract
\endtopmatter

\document
\head 0.  Introduction \endhead

 Let $L$ be a very ample line bundle on a smooth
 complex projective variety $X$, giving an embedding
 $X\subset\Bbb P^N$.  It is well-known that projective
 normality of the embedding (or {\sl normal
 generation} of $L$) is equivalent to the vanishing
 $$
 H^1 (\Bbb P^N ,\Cal I_X (k))=0, \qquad \text{all }k,
 $$
 where $\Cal I_X$ is the ideal sheaf defining $X$;
 further, all sufficiently high powers of $L$ are
 normally generated.  In this paper, we shall be
 concerned with the condition on $L$ (or its
 embedding)
 $$
   H^1(\Bbb P^N ,\Cal I_X^2  (k))=0,
 \qquad \text{all } k\neq 2. \tag"($\ast$)"
 $$
 (For k=2, Proposition 1.8 shows this group is
 frequently the kernel of the Gaussian map of $L$,
 hence is rarely 0.)  This question arises naturally
 because these cohomology groups give the torsion
 submodule of the K\"ahler differentials of the affine
 cone over $X$ (Proposition 1.4).  Our first main
 results are:

 \proclaim{Theorem 2.1}  {\rm ($\ast$)} holds for $X=\Bbb
 P^n$ and
 any very ample $L$.\endproclaim

 \proclaim{Corollary 3.3}  Given any $X$ and very
 ample $L$,
   all sufficiently high powers of $L$ satisfy
 {\rm ($\ast$)}.\endproclaim

 \proclaim{Corollary 5.7}  Let $C\subset\Bbb P^{g-1}$ be the
 canonical embedding of a general non-hyperelliptic
 curve of genus $g\geq 3$.  Then {\rm ($\ast$)} holds. \endproclaim

 Verification of ($\ast$) for $\Bbb P^n$ surprisingly
 turns out to be non-trivial; even the $n=1$ case
 requires some serious thought! One must deduce via
 representation theory of $G=SL(n+1)$ and the methods
 of [W4] the surjectivity of
 $$
 \Gamma (\Cal I (k)) \rightarrow \Gamma(\Cal I/\Cal
 I^2(k)),
 \qquad k\geq 3.
 $$
 The difficulty is that the second space is a
 reducible $G$-module, and only half of its
 irreducible constituents are obviously in the image.
   Corollary 3.3
 is deduced from Theorem 2.1.  In general, ($\ast$) is difficult to
 verify; it holds for the Pl\"ucker embedding of a Grassmannian
$G(2,n+1)$, and for an embedding of a curve via a line bundle of degree
 $\geq 2g+3$ (J. Rathmann \cite{R}, unpublished).
 \par The main point of Corollary 5.7 is to calculate
 for $C$ which is {\sl pentagonal}, i.e., with a
 base-point-free $g^1_5$ (Theorem 5.3).  Such a $C$
 sits naturally on a 4-dimensional rational scroll
 $X$, for which ($\ast$) must be proved.  It is then a
 delicate calculation to \lq\lq descend" this result to
   $C$; the key point is that ($\ast$) holds for 5 points in general
position in $\Bbb P^3$.
 As ($\ast$) fails for the canonical embedding of a
 tetragonal curve or plane sextic, one can
 optimistically hope to prove the following
 \proclaim{Conjecture}  Let $C$ be a curve of
   Clifford
 index $\geq 3$.  Then the canonical embedding of $C$
 satisfies {\rm ($\ast$)}.\endproclaim

 This paper is motivated by the question of whether
 a canonical curve $C\subset\Bbb P^{g-1}$ is {\sl
 extendable}, i.e., is a linear section of an
 $X\subset\Bbb P^g$ which is not the cone over $C$.
 Such an $X$ is {\sl canonically trivial} (c.t.):
 normal, $\omega_X\cong\Cal O_X$, and $h^1(\Cal O_X)=0$
 (\cite{E}, \cite{W5}).  If smooth, $X$ is a K-3;
 but any normal quartic in $\Bbb P^3$ is c.t. In \cite{W2},
  extendability of $C$ is studied via the deformation
 theory of the affine cone $A$ over $C$.
  Then $T^1_A$, the module of first-order deformations, is non-trivial
 in degree $-1$; by duality, this means that the Gaussian map
 $$
 \Phi_K: \wedge ^2 \Gamma (C,K) \rightarrow \Gamma
 (C,K^{\otimes 3})
 $$
 is not surjective. Lifting a first-order deformation
 of $A$ to a higher order requires control of the
 obstruction space
 $T^2_A$ which, by local duality, is dual to torsion
 in the K\"ahler differentials.  Noting the grading on
 these modules, and combining (1.6) with (6.4) yields the
 \proclaim{Theorem}  Let $C\subset\Bbb P^{g-1}$ be a
 canonical curve satisfying {\rm ($\ast$)}, and let $A$ be
 the affine cone.
 \roster
 \item"(a)" If $3\leq g \leq 10$, $g\neq 9$, then
 $T^2_A =0$.\endroster \vskip1ex \par \roster
 \item "(b)" If $g=9$ or $g\geq 11$, then $T^2_A$ is
 concentrated in degree $-1$.
 \endroster
 \endproclaim
 (The $g=9$ exception is the non-injectivity of $\Phi_K $,
 shown in [C-M] and reproved in (6.5.2)
 using work of S. Mukai). We use this to prove a
 partial converse to the non-surjectivity of the
 Gaussian for a K-3 curve:
 \proclaim{Theorem 7.1}  Let $C\subset\Bbb P^{g-1}$
 be a canonical curve satisfying {\rm ($\ast$)}.  Then $C$
 is extendable iff the Gaussian $\Phi_K $ is not
 surjective.\endproclaim
   The canonical embedding of a smooth plane curve $C$
 of degree $\geq 7$ satisfies ($\ast$), with corank
 $\Phi_K=10$, hence is extendable; we give an
 explicit description in (7.3), where a typical extension $X$ has a
 non-smoothable simple elliptic singularity.  Such a
 $C$ sits on no K-3 surface, by \cite{G-L} or
 \cite{W5}.  Because of the special geometry of c.t.
 surfaces, we offer the
 \proclaim{Conjecture} A Brill-Noether-Petri general
 curve of genus $\geq 8$ sits on a K-3 surface if and
 only if the Gaussian is not surjective. \endproclaim

 By (1.13), condition ($\ast$) is a natural one to
 consider for a line bundle already known to be {\sl
 normally presented} (i.e., satisfying M. Green's
 condition $(N_1)$ --- projectively normal, with
 homogeneous ideal generated by quadrics).
 Cubics vanishing twice on $C$ vanish on each secant line, hence on the
 secant variety Sec$(C)$; so ($\ast$) should be related to Sec$(C)$
being defined by cubics.
  Following work of
 Aaron Bertram \cite{B} and Michael Thaddeus \cite{T},
 one should study more generally for projectively
 normal $C\subset \Bbb P^N$ the spaces
 $$
 H^0(\Bbb P^N,\Cal I^n_C(n+m)).
 $$

  \vskip1ex
 \par The paper is organized as follows: Section 1
 introduces the groups $H^1(\Cal I^2(k))$, relating
 them to cones and to Gaussian mappings; and ($\ast$)
 is compared with other nice properties of high powers
 of ample line bundles, specifically normal
 presentation.  In Section 2, ($\ast$) is proved for
 all embeddings of $\Bbb P^n$ and discussed for
 complex homogeneous spaces.  How to \lq\lq descend" the
 property ($\ast$) from a variety to a subvariety is
 discussed in \S3; products and scrolls are considered in \S4,
 allowing one to prove ($\ast$) for linear sections of
 a Segre embedding.  Calculations for pentagonal and
 tetragonal curves are done in Section 5, using
 results about the scrolls on which they naturally sit
 (as described in [S]).  What does or might occur for
 other canonical curves is discussed in \S6, where the
 work of Mukai on curves with $g\leq 9$ gives another
 approach to the condition ($\ast$).
  Section 7 describes the application to extendability
 of canonical curves.
 \par This paper has benefitted from conversations with Rob Lazarsfeld
and Michel Brion, as well as from careful comments from the referee.
Research was partially supported by
 DMS-9302717.

 \head 1.  A cohomological vanishing
 condition\endhead

 (1.1)  Let $X\subset\Bbb P^N$ be a (non-degenerate)
 projectively normal subvariety, with ideal sheaf
 $\Cal I=\Cal I_X$; thus, $H^1(\Cal I(k))=0$, all $k$.
 There is an exact sequence
 $$
 \Gamma (\Cal I(k))\rightarrow \Gamma (\Cal I/\Cal
   I^2(k))\rightarrow H^1(\Cal I^2(k))\rightarrow 0.
 \tag1.1.1
 $$
 We study the additional condition
 $$
 H^1(\Cal I^2(k))=0, \qquad k\geq 3.  \tag"($\ast$)"
 $$
 By (1.3.2), $H^1(\Cal I^2(k))=0$ is nearly automatic
 for $k\leq 1$, but unusual for $k=2$.  One has the easy
 \proclaim{Proposition 1.2} If $X\subset\Bbb P^N$  is
 a complete intersection, then
 $$
 H^1(\Cal I^2(k))=0, \qquad \text{all } k.
 $$
 \endproclaim

 \flushpar {\it Proof.}  If the degrees of the defining
 equations are $\{d_i\}$, one has a surjection
 $$
 \oplus\Cal O_{\Bbb P}(-d_i) \rightarrow \Cal I
 \rightarrow 0;
   $$
 tensoring with $\Cal O_X$ induces an isomorphism
 $$
 \oplus  \Cal O_X(-d_i) \cong \Cal I/\Cal I^2.
 $$
 Since $X$ is projectively normal, $\Gamma
 (\Cal O_{\Bbb P} (k-d_i))\rightarrow\Gamma(\Cal O_X
 (k-d_i))$ is surjective for all $i$.  Now use (1.1.1).

\vskip1ex

 (1.3)  If $X\subset\Bbb P^N$  is non-singular and
 linearly normal, $L=\Cal O_X (1)$, and $V=\Gamma(\Cal O_{\Bbb P}(1))=
 \Gamma(\Cal O_X(1))$,
there is a basic diagram
 $$\matrix\format \c&&\quad\c \\
 && 0  && 0\\
   \vspace{2\jot}
 &&\downarrow  &&\downarrow\\
   \vspace{2\jot}
 && \Cal I/\Cal I^2  &\cong  &\Cal I/\Cal I^2\\
   \vspace{2\jot}
 &&\downarrow  &&\downarrow\\
   \vspace{2\jot}
 0  &\rightarrow  &\Omega_{\Bbb P}^1 \otimes \Cal O_X
 &\rightarrow  &V\otimes \Cal O_X(-1)
 &\rightarrow &\Cal O_X  &\rightarrow  &0\\
   \vspace{2\jot}
 &&\downarrow  &&\downarrow &&\Vert \\
   \vspace{2\jot}
 0  &\rightarrow  &\Omega_X^1  &\rightarrow  &\Cal E
 &\rightarrow  &\Cal O_X &\rightarrow  &0\\
 \vspace{2\jot}
 &&\downarrow  &&\downarrow\\
   \vspace{2\jot}
 &&0  &&0
 \endmatrix\tag1.3.1
 $$
 The middle row is the Euler sequence of $\Bbb P^N$
 restricted to $X$; it may be described intrinsically
 as a twist of
 $$
 0 \rightarrow  M_L \rightarrow \Gamma (L)\otimes
 \Cal O_X \rightarrow L \rightarrow 0,
 $$
 where the second map is evaluation; thus,
 $M_L=\Omega_{\Bbb P}^1 (1)\otimes \Cal O_X$ (cf.
 [L]).  The bottom row is the extension corresponding
 to the class of the line bundle $L$ in
 $H^1(\Omega^1_X)$.  It is clear that $\Gamma
 (\Omega_{\Bbb P}^1 (k)\otimes \Cal O_X ) = \Gamma
 (M_L\otimes L^{k-1})=0$, $k\leq 1$, whence $\Gamma
 (\Cal I/\Cal I^2(k))=0$ in that range.  From (1.1.1),
 $$\gather
 X\subset\Bbb P^N\ \text{linearly normal and
 nonsingular
  implies} \tag1.3.2\\ H^1(\Cal I^2(k))=0,\qquad\text{all }
  k\leq 1.
 \endgather$$

 \proclaim{Proposition 1.4}  Let $X\subset\Bbb P^N$ be a
   non-singular projectively normal variety, with affine
   cone $A=\oplus  \Gamma(X,\Cal O_X (n))$.  Let
 $\Omega^1_A
   = \oplus (\Omega^1_A)_{{}_k}$ be the graded module of K\"ahler
 differentials of $A$.  Then for all $k$,
 $$
 (\text{{\rm Tors }}\Omega^1_A)_{{}_k} \cong H^1(\Cal
 I^2(k)). \tag 1.4.1
   $$
 \endproclaim
 \flushpar {\it Proof.} Write $P=\oplus\Gamma(\Cal O_{\Bbb P}(k)),
 I=\oplus I_k=\oplus\Gamma(\Cal I(k))$;
  then $A=P/I=\oplus\Gamma(\Cal O_X(k))$. The K\"ahler differentials of $A$
satisfy
 $$
 I/I^2 \rightarrow \Omega^1_P \otimes A
 \rightarrow \Omega^1_A \rightarrow 0;
 $$
 the first map is induced from exterior
  differentiation $d:I\rightarrow\Omega^1_P$, and its
   kernel is a torsion A-module (A has an isolated singularity). The
graded module $\Omega^1_P$ is canonically $V\otimes_{\Bbb C}P(-1)$, so
$$
 (\Omega^1_P \otimes A)_k \simeq V\otimes_{\Bbb C} A_{k-1}.
$$
 The kth graded piece of $d(I/I^2)$ is the image
   of the composition
 $$
 I_k =\Gamma (\Cal I(k)) \rightarrow \Gamma (\Cal
   I/\Cal I^2(k))\subset \Gamma (\Omega^1_{\Bbb P}
   \otimes \Cal O_X (k))\subset
 V\otimes A_{k-1}.
   $$
 Letting $J_k$ denote the image of $I_k$ in $\Gamma
   (\Cal I/\Cal
 I^2(k))$, we thus have
 $$\gather
 (\Omega^1_A)_{{}_k} = (\Omega^1_P \otimes A)_k /J_k\\
 0 \rightarrow J_k \rightarrow \Gamma (\Cal I/\Cal
   I^2(k)) \rightarrow H^1(\Cal I^2(k)) \rightarrow 0.
   \tag1.4.2\endgather$$
 \par Spec $A$ is obtained by collapsing the
   0-section $X$ of the geometric line bundle associated to $L^{-1},\ \pi
   :V \rightarrow X$.  Write $U=V-X=\text{
   Spec}(A)-\{0\}$.  Since $\pi^\ast \Cal E \cong
   \Omega^1_V (\log X)$, (e.g., \cite{W1}, 3.3), one has
 $$
   \Gamma(U,\Omega^1_A )
   =\Gamma(U,\Omega^1_V)
   =\oplus\Gamma(X,\Cal E(k)), \qquad \text{all
   }k\in\Bbb Z .
   $$
  The local cohomology sequence for $U\subset \text{Spec}\,A$ gives
$$
(\text{Tors }\Omega^1_A)_k = \text{Ker}\{ V\otimes\Gamma(\Cal
O_X(k-1))/J_k
 \rightarrow \Gamma(\Cal E(k))\}.
$$
Thus, the inclusion
$$
V\otimes \Gamma(\Cal O_X(k-1))/\Gamma(\Cal I/\Cal I^2(k))
\hookrightarrow
\Gamma(\Cal E(k))
$$
given by the middle column of (1.3.1), plus the sequence (1.4.2),
produces (1.4.1).
 \vskip1ex
  \flushpar {\it Remark} (1.5.1). Proposition (1.4)
 remains true for
 singular $X$ if $\Omega^1_X$ is reflexive and $\Cal
 I/\Cal I^2$ is
  torsion-free (e.g., $X$ is a normal
   local complete intersection).
  \vskip1ex \flushpar  (1.5.2) Higher local cohomology of
 $\Omega^1_A$ is computed
  more easily:
 $$
   H^i_{\{0\}}  (\Omega^1_A)_{{}_k} \cong
   \biggl\{
   \aligned &\text{Coker } \Phi(L, L^{k-1} ):
   \Gamma (\Omega^1_{\Bbb P} \otimes \Cal O_X
 (k))\rightarrow
   \Gamma (\Omega^1_X  (k))  \qquad i=1\\
   &H^{i-1} (X,\Cal E(k)) \qquad i\geq 2.\endaligned
  $$
 For the first assertion, see [W3],(1.3.5),
   plus (1.7) below for the definition of the Gaussian
 $\Phi$;
   the second follows from computing $H^{i-1}
 (U,\Omega^1_V)$ as above.
 This yields the local cohomology
   of the $A$-module $I/I^2$, except for its torsion.
   Since the kernel of $I_k \rightarrow J_k$ is
 $\Gamma(\Cal I^2(k))$,
 $$\gather
   \text{{\rm Tors }}(I/I^2)_k=\Gamma (\Cal
 I^2(k))/(I^2)_k,\\
 \intertext{whence}
 \text{{\rm Tors }}(I/I^2)_3=\Gamma (\Cal I^2(3)).
   \endgather
   $$
 \proclaim{Corollary 1.6}  Let $C\subset\Bbb P^{g-1}$ be the
 canonical embedding of a non-hyperelliptic curve of
 genus $g\geq 3$, with $A$ the coordinate ring of the
 affine cone.  Then the graded pieces of the cotangent
 modules $T^1_A$ and $T^2_A$ are given by
 $$\align
   (T^1_A)_{{}_k} &\cong \text{{\rm Coker }}
 \Phi(K,K^{-k})^\ast\tag1.6.1\\ \vspace{1\jot}
 (T^2_A )_{{}_k} &\cong H^1(\Cal I^2(1-k))^\ast. \tag1.6.2
   \endalign$$
 In particular, if $H^1(\Cal I^2(k))=0$, $k\neq 2$,
   then $T^2$ is concentrated in degree $-1$, and is
 dual to the kernel
   of the Gaussian  $\Phi_K :\wedge^2\Gamma(K)
 \rightarrow\Gamma(K^{\otimes 3})$.
   \endproclaim
 \flushpar {\it Proof.} The cone is a Gorenstein surface
 singularity, with dualizing module $\omega_A=A(1)$.
   By local duality, one has ([W2], (3.8))
 $$\align
 (T^1_{1 - k} )^\ast &\cong H^1_{\{ 0\}}
 (\Omega^1_A)_{{}_k}\\ \vspace{1\jot}
 (T^2_{1 - k} )^\ast &\cong H^0_{\{ 0\}}
 (\Omega^1_A)_{{}_k}.
   \endalign$$
 Now use Proposition 1.4 above and Proposition 1.8
 below.
\vskip1ex

  \par (1.7).
  Let $L$ and $L'$ be two line bundles on an arbitrary
   smooth projective $X$. The {\sl diagonal
 filtration} $\Cal R_i (L,L')$
   of $\Gamma (L)\otimes\Gamma (L')$ is defined by
   $$
   \Cal R_i(L,L')\equiv\Gamma(X\times X,\Cal
 I^i_{\triangle}
   (L\boxtimes L'));
   $$
 where $\Cal I_\triangle$ is the ideal sheaf of the
 diagonal on $X\times X$.
 Since $\Cal I_\triangle^i/\Cal I_{\triangle}^{i+1}
 \cong S^i\Omega^1_X$,
   $\Cal R_{i+1}$ is the kernel of a {\sl Gaussian map}
 $$
   \Phi_i(L,L'):\Cal
 R_i(L,L')\rightarrow\Gamma(X,S^i\Omega^1_X
 \otimes L\otimes L'),
   $$
 given by restriction to the diagonal on $X\times X$.
 Writing multiplication as $\mu\equiv\Phi_0$, and
 $\Phi\equiv\Phi_1$, one has
 $$
   \Cal R (L,L') \equiv \Cal R_1(L,L')=\text{ Ker }(\Phi
 (L,L'):
 \Gamma (L)\otimes \Gamma (L')\rightarrow\Gamma
 (L\otimes L')).
   $$
 If $L$ is very ample (or more generally, generated by
 global sections), one has
   $$
   \Cal R(L,L')=\Gamma (X,M_L \otimes L'),
  $$
 and the first Gaussian map is given by
 $$
 \Gamma (X,M_L \otimes L')\rightarrow \Gamma (X,\Omega^1_X
 \otimes L\otimes L').
   $$
 One concludes:
 $$\gather
  0\rightarrow \text{ Coker } \mu (L, L^{k-1}
 )\rightarrow H^1 (M_L\otimes L^{k-1})\rightarrow
 \Gamma (L)\otimes H^1(L^{k-1}) \tag1.7.1 \\
   \Gamma (\Cal I/\Cal I^2(k))=\text{ Ker } \Phi(L,
 L^{k-1} )=\Cal R_2(L, L^{k-1} )\tag1.7.2\\
 0\rightarrow \text{ Coker } \Phi(L, L^{k-1}
 )\rightarrow H^1(\Cal I/\Cal I^2(k))\rightarrow
 H^1(M_L\otimes L^{k-1}).\tag1.7.3
  \endgather$$
 When $L=L'$, the diagonal filtration is stable by the natural involution
of \newline $\Gamma
 (L)\otimes\Gamma(L)$, so induces a filtration of both
 $S^2\Gamma (L)$ and $\wedge^2\Gamma (L)$ (viewed as subspaces); we set
 $$
 \Cal R_i^s=S^2 \Gamma (L)\cap \Cal R_i (L,L),\  \Cal R_i^a
 =\wedge^2 \Gamma (L) \cap \Cal R_i (L,L).
   $$
 The Gaussians $\Phi_i$ are $(-1)^i$-symmetric, hence for odd $i$ vanish on
 symmetric tensors, and for even $i$
 vanish on the alternating ones; therefore,
   $$
   \Cal R_{2i}^s =\Cal R_{2i-1}^s \ \
   \text{and}\ \  \Cal R^a_{2i} =\Cal R^a_{2i+1}.
   $$
 In particular,
 $$\align
 \Cal R(L,L) &=\Cal R_1^s\oplus\Cal R_1^a = I_2 \oplus
 \wedge^2 \Gamma (L)\\
  \Cal R_2(L,L) &=\Cal R ^s_2 \oplus \Cal R ^a_2 =I_2
 \oplus \text{Ker }\Phi_L ,
   \endalign
   $$
 where $\Phi_L:\wedge^2 \Gamma(L)\rightarrow \Gamma
 (\Omega^1 \otimes L^2)$ is the usual Gaussian map
 $(s\wedge t\mapsto sdt-tds)$. Comparing the last
 decomposition with (1.7.2) for k=2, one identifies
 the cokernel of $I_2 \rightarrow \Gamma(\Cal I/\Cal
 I^2(2))$ with Ker $(\Phi_L )$.  Thus,
  \proclaim{Proposition 1.8}  Suppose $X\subset\Bbb
 P^N$ is a
   nonsingular projectively normal subvariety, with
 $L=\Cal O_X (1)$.  Then
 $$
   H^1 (\Cal I_X^2 (2))\cong \
   \text{{\rm Ker}}\  (\Phi_L :\wedge^2\Gamma
 (L)\rightarrow
   \Gamma (\Omega^1_X \otimes L^2 )).
   $$
   \endproclaim

  \flushpar {\it Remark} (1.9). It is thus clear from a dimension
 count using Riemann-Roch that for $L$ sufficiently
 ample on a given $X$,
   $H^1(\Cal I^2(2))\neq 0$.

   (1.10) Recall there are natural functorial maps
 $$
   \Cal R_i (L,L')\otimes \Gamma (L'')\rightarrow
   \Cal R_i(L,L'\otimes L''),
   $$
 which commute with the Gaussians.  We consider for a
 very ample $L$ the surjectivity of the induced maps
 $$
   \Cal R_i^p \otimes \Gamma (L^{k}) \rightarrow
   \Cal R_i (L,L^{k+1}),\tag1.10.1
   $$
 where $p=s$ (respectively $p=a$) if $p$ is even (resp. $p$ is odd). Since
   $\Cal R^s_0=S^2\Gamma(L)$ and $\Cal
 R^a_1=\wedge^2\Gamma(L)$, the first
   two cases of (1.10.1) are covered by
 \proclaim{Proposition 1.11} Let $L$ be a very ample
 line bundle on $X$.  Then
   \roster \item"(a)" $L$ is normally generated
 iff\endroster
 $\qquad\qquad S^2\Gamma (L)\otimes \Gamma (L^{
 k}) \rightarrow \Gamma (L)\otimes \Gamma (L^{
 k+1})$ is surjective, $k\geq 1$.
   \roster \item"(b)" If $L$ is normally generated,
 then it is
 normally presented iff\endroster
 $\qquad \qquad\wedge^2\Gamma (L)\otimes\Gamma
 (L^{k}) \rightarrow \Cal R_1(L,L^{
 k+1})$ is surjective, $k\geq 1$.\endproclaim
 \flushpar {\it Proof.} The map in (a) factors through
   $\Gamma (L)\otimes \Gamma (L)\otimes \Gamma
 (L^k)$, so the surjectivity of the map in (a)
 easily implies normal generation.
 For the converse, one reduces by restriction to the
 case $X=\Bbb P^n,\ L=\Cal O(1)$; and, the result for
 $k=1$ will imply that for general $k$. One must
 therefore prove for a vector space $V$ the
 surjectivity of the composition
   $$
   S^2V\otimes V\subset V\otimes V\otimes
 V\twoheadrightarrow
   V\otimes S^2V. \tag1.11.1
  $$
 But an element of the kernel is symmetric with
 respect to the first two entries, anti-symmetric
   with respect to the last two.  The signature being the only
non-trivial character of $S_3$,
 there could be no such
 vectors.  The map in (1.11.1) is injective and thus surjective.
 \par The proof for (b) is slightly easier. For
 $X\subset\Bbb P,\ L=\Cal O_X (1)$, one examines the
 diagram
 $$\CD
  \wedge^2\Gamma (\Cal O_{\Bbb P}(1))\otimes \Gamma
 (\Cal O_{\Bbb P}(k))  @>>>  \Gamma (\Omega^1_{\Bbb P}
 (k+2))\\
 @VVV  @VVV  \\
 \wedge^2\Gamma (\Cal O_X (1))\otimes \Gamma (\Cal
 O_X(k))
 @>>>  \Gamma (\Omega^1_{\Bbb P} \otimes \Cal O_X
 (k+2)). \endCD$$
 By (1.3) and (1.6), the horizontal maps are the maps
 $\Cal R^a_1(L)\otimes\Gamma(L^{k})
 \rightarrow \Cal R_1(L,L^{k})$. The top map
 is surjective for $k\geq 1$ (e.g. it is a non-0
 equivariant map into an irreducible G-bundle, where
 $G=\, \text{Aut}\,(\Bbb P))$. The left map is clearly
 surjective.
 The cokernel of the right vertical map is (by
 tensoring the Euler sequence)
 $$
 H^1(\Bbb P,\Omega^1_{\Bbb P} (k+2)\otimes \Cal
 I_X)\cong \, \text{Coker}\, \{\Gamma(\Cal
 I_X(k+1))\otimes \Gamma (\Cal O_{\Bbb
 P}(1))\rightarrow \Gamma (\Cal I_X (k+2)\}.
 $$
 This cokernel vanishes for $k\geq 1$ iff
 $X\subset\Bbb P$ is ideal-theoretically generated by
 quadrics.
\vskip1ex

 \flushpar {\it Remarks} (1.12.1) For $i=0$ or 1, the
 maps in (1.10.1)
 are surjective iff they are surjective without the
 \lq\lq p".  For $i=1$, this condition for normal
 presentation is already found in [Mu].
 \par (1.12.2)  The previous result for $i=0$ or 1
 motivates the following
 \proclaim{Proposition 1.13}  Suppose $X\subset\Bbb
 P^N$ is smooth and normally presented.  Then for
 $k\geq 1$,
 $$
 H^1(\Cal I^2(k+2)) \cong\,\text{{\rm Coker }} (\Cal
 R^s_2 \otimes \Gamma (L^k) \rightarrow \Cal
 R_2(L,L^{k+1})).
 $$
 \endproclaim
 \flushpar {\it Proof.} Use (1.1.1), (1.7.2), and the maps
 $$
 I_2 \otimes \Gamma (L^k) \twoheadrightarrow I_{k +2}
 \rightarrow \Gamma (\Cal I/\Cal I^2(k+2)).
 $$

 (1.14)  We state without proof an unpublished
 theorem of J\"urgen Rathmann which allows calculation
 of $H^1(\Cal I^2(k))$ in the case of a curve $C$.  On
 $C\times C\times C$, let $\Cal I_{12},\ \Cal I_{13},\
 \Cal I_{23}$ be the ideal sheaves of the \lq\lq
 partial diagonals", i.e., the divisors for which the
 corresponding coordinates are equal; also, $\pi_1: C\times C\times C\to
 C$ is projection.

 \proclaim{Theorem 1.15 (Rathmann [R])} Let
 $C\subset\Bbb P^N$ be a projectively normal embedding
 of a smooth projective curve, ideal-theoretically
 defined by quadrics.  Write $L=\Cal O_C(1)$.  Then
 there is an
 exact sequence on $C$
 $$
  0\rightarrow \wedge^2 M_L\otimes L
  \rightarrow \pi_{1\ast} (\Cal I_{12}\Cal I_{23}\Cal
 I_{13}
  (L\boxtimes L \boxtimes L)) \rightarrow I_2 \otimes
 L\rightarrow \Cal I/\Cal
  I^2(3)
  \rightarrow 0.
  $$\endproclaim
   \proclaim{Corollary 1.16}  Suppose {\rm deg} $L\geq
 2g+3$.  Then the
 embedding defined by $L$ satisfies $H^1(\Cal
 I^2(k))=0$, all $k\geq
  3$.
 \endproclaim

 \flushpar {\it Proof}. Use (1.13), (1.15), and a vanishing result
$$
H^1 (C\times C\times C,\ I_{12}I_{23}I_{13}
 (L_1 \boxtimes L_2\boxtimes L_3))=0,
$$
if deg$\,L_i\geq 2g+3$, proved along the general lines of \cite{L}.

  \vskip3ex
  \head 2.  Projective space and other G/P \endhead

   \proclaim{Theorem 2.1}  Let $X=\Bbb P^n\subset\Bbb
 P^N$ be the Veronese embedding of projective space
 defined by $L=\Cal O_{\Bbb P^n}(r)$, where $r\geq 1$.
 Then
 $$
 H^1(\Bbb P^N , \Cal I^2(k))=0, \qquad \text{for
  }k\neq 2.
 $$ \endproclaim

 \flushpar {\it Proof.} Bypass the trivial case
 $r=1$.  It is
  standard that these embeddings are projectively
 normal and normally presented.  By (1.3.2), one may
 assume $k\geq 3$.  By (1.13), it suffices to prove
 $$
 \Cal R^s_2\otimes\Gamma (L^k) \twoheadrightarrow \Cal
 R_2(L,L^{k+1}), \qquad\text{all } k\geq 1.
 $$
 In fact, we will prove by a descending induction on $i$ that
 $$
 \Cal R^s_{2i}\otimes \Gamma (L^k) \twoheadrightarrow
 \Cal R_{2 i} (L,L^{k+1}), \qquad\text{all }i,k\geq
 1.\tag2.2.1
 $$
 \par (2.3) $\Bbb P^n$ is a homogeneous space for
 $G=SL(n+1)$; $G$ acts on all relevant spaces, and all
 natural maps are $G$-equivariant.  One deduces, from
 known results or from ([W4], \S3,4):
 \vskip1ex \flushpar (2.3.1) $S^i \Omega^1_{\Bbb P}(j)=\pi_\ast
 L_{i,j}$, where $L_{i,j}$ is a line bundle on the
 flag manifold $G/B$, \roster\item"" and $\pi$ is the projection $G/B
\rightarrow\Bbb P^n$.
\endroster

 \vskip1ex
 \flushpar (2.3.2)  $\Gamma (S^i\Omega^1_{\Bbb P}(j))$ is
 an irreducible $G$-module, non-0 iff $j\geq 2i$.

 \vskip1ex
 \flushpar (2.3.3)  For $j\geq 2i,\ j'\geq 2i'$, one has a
 surjection
 $$
 \Gamma (S^i\Omega^1_{\Bbb P}(j))\otimes\Gamma
 (S^{i'}\Omega^1_{\Bbb P}(j')\rightarrow\Gamma(S^i\Omega^1_{\Bbb
 P}\otimes S^{i'}\Omega^1_{\Bbb P}(j+j'))\rightarrow\Gamma
 (S^{i+i'}\Omega^1_{\Bbb P} (j+j')).
 $$

 \flushpar (2.3.4) $\Cal R_i (L,L^k)=0$ iff $i>r$  (recall
 $L=\Cal O(r))$.
 \vskip1ex
 \flushpar (2.3.5) For all $i\leq r,\ \Phi_i:\Cal R_i (L,L^k) /\Cal R_{i+1}
 (L,L^k) \overset{\sim}\to{\rightarrow} \Gamma (S^i\Omega^1\otimes
 L^{k+1})$ \vskip1ex
  \flushpar is an isomorphism of irreducible $G$-modules.
  In particular,
   $$\gather
   \Cal R^s_2 = I_2 = \underset {i \geq 1}\to{\oplus}
 \Gamma (S^{2i}
  \Omega^1 \otimes L^2 )\\
   \Cal R_2 (L,L^{k+1} )= \underset {i \geq
 2}\to{\oplus} \Gamma
  (S^i\Omega^1 \otimes L^{k+2} ).
  \endgather$$
 While it is clear that
   $$
   \Gamma (S^{2i} \Omega^1 \otimes L^2 )\otimes \Gamma
  (L^k) \twoheadrightarrow \Gamma (S^{2i}\Omega^1
  \otimes L^{k+2})
   $$
   is surjective for $2i\leq r$ (2.3.2), the heart of
  the theorem is to prove that the \lq\lq odd part" of
  $\Cal R_2(L,L^{k+1})$ is also in the image.

  (2.4) We consider the diagram
 $$\matrix\format \c&&\quad\c \vspace{3\jot}\\
   0  &\rightarrow  &\Cal R^s_{2i+1}\otimes \Gamma (L^k)
 &\rightarrow  &\Cal R^s_{2i}\otimes \Gamma (L^k)
 &\rightarrow  &\Gamma (S^{2i}\Omega ^1\otimes
 L^2)\otimes \Gamma (L^k)  &\rightarrow  &0\\
 \vspace{3\jot}
   &&\sideset\and{\alpha_i}\to{\downarrow}
 &&\sideset\and{\beta_i}\to{\downarrow}
 &&\sideset\and{\gamma_i}\to{\downarrow}\\
 \vspace{3\jot}
   0  &\rightarrow  &\Cal R_{2i+1} (L,L^{k+1})
 &\rightarrow  &\Cal R_{2i} (L,L^{k+1})
 &\rightarrow  &\Gamma (S^{2i}\Omega ^1\otimes
 L^{k+1})  &\rightarrow  &0.
  \endmatrix
  $$
 The top row is exact by (2.3.5) and the equality
  $\Cal R^a_{2i+1}=\Cal R^a_{2i}$; the bottom, by
  (2.3.5). $\gamma_i$ is
   surjective by (2.3.3).  If $2i\geq r$, then $\Cal
  R_{2i+1}=0$, so $\beta_i$ is surjective.

  \par Now assume
  that $2i <r$ and $\beta_{i+1}$ is
surjective; we prove $\beta_i$ is surjective. By the snake lemma, it
suffices to check the
 \proclaim{Claim (2.5)} {\rm
  Ker}$(\gamma_i)\rightarrow$ {\rm Coker} $(\alpha_i)$
  is surjective for $0 <2i<r$.
  \endproclaim
 \par (2.6) First, since
  $\Cal R^s_{2i+1} =\Cal R^s_{2i+ 2}$, Im$(\alpha _i)
  =$ Im$(\beta_{i+1})=\Cal R_{2i+2} (L,L^{k+1})$ (by
  induction), one has
   $$
 \text{{\rm Coker}} \ (\alpha _i)\cong \Cal R_{2i+1}
  (L,L^{k+1}) /\Cal R_{2i+2} (L,L^{k+1}) \cong
   \Gamma (S^{2i+1}\Omega^1\otimes L^{k+2} ).
   $$
  To identify the map in (2.5), we
  first define bundles $M_k$ by
 $$\gather
 0\rightarrow M_k\rightarrow \Gamma (L^k)\otimes
  \Cal O_X \rightarrow L^k\rightarrow 0.\\
   \intertext{Then}
   \text{Ker }(\gamma_i) =\Gamma (S^{2i} \Omega^1
  \otimes L^2 \otimes M_k).
  \endgather$$
 Via the surjection $M_k\twoheadrightarrow \Omega^1_X
  \otimes L^k$ (1.3.1), the map in (2.5) is the composition
 $$
   \Gamma(S^{2i}\Omega^1\otimes L^2\otimes M_k)
 \rightarrow \Gamma (S^{2i}\Omega^1\otimes L^2
 \otimes\Omega^1\otimes L^k)\rightarrow \Gamma
 (S^{2i+1}\Omega^1\otimes L^{k+2} ).
   $$
 We assert that one has a surjection after composing
 this map at the beginning with
   $$
   \Gamma(S^{2i}\Omega^1\otimes L^2(-1))\otimes
  \Gamma(M_k(1))\rightarrow\Gamma (S^{2i}\Omega^1\otimes L^2
  \otimes M_k).
   $$
 To accomplish this, note first that
  $$
 \Gamma(M_k(1))\rightarrow\Gamma(\Omega^1\otimes L^k(1))
  $$
 is the Gaussian on $\Bbb P^n$ for $\Phi (L^k,\Cal
 O_{\Bbb P}(1))$, hence is surjective.  Next, $2i<r$
 implies $4i<2r-1$; by (2.3.2),
   $$
   \Gamma (S^{2i} \Omega^1 \otimes L^2 (-1))\neq 0.
   $$

 Finally,
   $$
   \Gamma (S^{2i} \Omega^1 \otimes L^2 (-1))\otimes
   \Gamma (\Omega^1\otimes L^k(1))\rightarrow \Gamma
  (S^{2i+1} \Omega^1 \otimes L^{k+2} )
   $$
   is surjective by (2.3.3).  This completes the proof
  of the claim, and hence of Theorem 2.1.

  \subhead Problem 2.7 \endsubhead  Let $L$ be a very
 ample line bundle on a complex homogeneous space
 $G/P$.  Does the corresponding projective embedding
 have the property of Theorem 2.1?
 \vskip1ex
  \par (2.8)  This seems a difficult question in
  general. It is a theorem of B. Kostant that the
  embedding defined by $L$ is normally presented, and
  (3.3) below gives the desired property for $L^k$,
  where $k$ is large.
 There is one class of \lq\lq easy" cases.

   \proclaim{Proposition 2.9}  Let $X=G/P\subset\Bbb
  P^N$ be an  embedding of a complex homogeneous space
  for which $\Cal I/\Cal I^2$ is an irreducible
  $G$-bundle.  Then
 $$
 H^1 (\Cal I^2_X  (k))=0, \qquad \text{all } k.
 $$ \endproclaim
 \flushpar {\it Proof.}  For $k\geq 2$, the non-0
 $G$-linear map $\Gamma (\Cal I(k))\rightarrow \Gamma
 (\Cal I/\Cal I^2(k))$ must be surjective, as the
 target space is an irreducible $G$-module.
 \vskip1ex
   \par (2.10)  Recall ([W4], \S4) that if
 $\lambda $ is a dominant weight for $G$,
 $L=L_\lambda$ the corresponding ample line bundle on
 $G/P$, there is a natural $P$-filtration $\{F_k\}$ of
 $\Gamma (G/P,L_\lambda)$, of length $\ell(\lambda)$.
 Now, $\Cal I/\Cal I^2$ as a G-bundle corresponds to
 the $P$-module $F_2$, the second piece of the
 filtration; its irreducibility therefore implies that
 $\ell(\lambda)=2$. (We exclude the case
 $\ell(\lambda)=1$, which occurs only for $(\Bbb
  P^n,\Cal O(1))$.)  It would  appear that the
 following is a complete list of weights with
 $\ell(\lambda)=2$ (cf. [W4], 4.7, although the
 entries for $B_n$ and $C_n$ were there incorrectly
 interchanged):
   $$\align
   A_n\negthickspace :\  &\omega_2,\  \omega_{n-1},\
  \omega_1+\omega_n,\  2\omega_1,\  2\omega_n\\
 B_n\negthickspace :\  &\omega_1\\
 D_n\negthickspace :\  &\omega_1;\ \omega_3,\
 \omega_4 \ (n=4);\
  \omega_4,\ \omega_5\  (n=5)\\
 E_6\negthickspace :\  &\omega_1,\ \omega_6
  \endalign
  $$
   Now, for the adjoint representation
  $\omega_1+\omega_n$ of $A_n$, the space $F_2$ is not
  irreducible, although (2.9) still holds (via (4.5)
  below).  One can check in every other case that
 $F_2$
  is indeed an irreducible $P$-module;  for instance,
  $\omega_1$ for $B_n$ and $D_n$ correspond to quadric
  hypersurfaces.  We record carefully one other case.

 \proclaim{Theorem 2.11}  Let $X=G(2,n+1)\subset\Bbb
  P^N$ be the Pl\"ucker embedding of the Grassmannian
 of lines in
  $\Bbb P^n$.  Then $H^1(\Cal I_X^2(k))=0$, all $k$.
  \endproclaim
   \flushpar {\it Proof.}  We show that $\Cal I/\Cal
 I^2$ is an
  irreducible $G$-bundle.  One considers the
 $G$-module
  $\Gamma (L_\lambda )= \wedge^2 (\Bbb C^{n+1})$ as a
  $P$-module, where $P$ is the parabolic which
  stabilizes a particular 2-dimensional subspace. As a
  module over a Levi component $L\cong SL(2)\times
  SL(n-1)\subset P$, it decomposes as a sum of three
  $L$-irreducibles
   $$
 \wedge^2 (\Bbb C^2\oplus  \Bbb C^{n-1} )\simeq \wedge^2
 \Bbb C^2 \ \oplus\ \Bbb C^2\otimes \Bbb C^{n-1}\ \oplus\
 \wedge^2 \Bbb C^{n-1}.
   $$
   The third term must be $F_2$, so it is irreducible
 as a $P$-module.

 \subhead Example 2.12\endsubhead  Let
 $X^6=G(2,5)\subset\Bbb P^9$ be the Grassmannian of
 lines in $\Bbb P^4$, which is defined by the
  Pfaffians of a generic skew-symmetric $5\times 5$ matrix.
  For the affine cone $A=P/I$, $I/I^2$ is a
  Cohen-Macaulay A-module, as can be proved by the methods of this
 paper; but J. Herzog \cite{H} has proved $I/I^2$ is always
Cohen-Macaulay when $A$ is Gorenstein of codimension 3. That $\Gamma
(\Cal I^2(3))$ is 0 (1.5.2) is also seen because such cubics vanish on
the secant lines of $X$, which span $\Bbb P^9$.

\vskip1ex
\par (2.13)  An unusual example in commutative algebra is also furnished
from (2.10).

 \proclaim{Theorem 2.14} Let $X\subset \Bbb P^{15}$ be
 the 10-dimensional homogeneous space
 corresponding to a half-spin
 representation $\omega_5$ of {\rm Spin(10)} $=D_5$.  Let $A=P/I$ be the
homogeneous coordinate ring, which is Gorenstein of
 codimension 5.  Then
 \roster
 \item "(a)"  $I/I^2$ is a Cohen-Macaulay $A$-module.\endroster

\vskip1ex \roster
 \item "(b)"  $\Omega^1_A$ is a Cohen-Macaulay
 $A$-module.
 \endroster
 \endproclaim

 \flushpar {\it Proof}. (Outline)  The half-spin
 representation $\Gamma (\Cal O_X(1))$ is a Spin(10)
 representation arising from a 10-dimensional complex
 vector space with a non-degenerate symmetric form.
 If $V$ is a maximal isotropic subspace, the
 stabilizer has Levi component SL(5), and the
 representation space decomposes as
 $$
 \Lambda ^0 V\oplus \Lambda ^2V \oplus \Lambda ^4V.
 $$
 The conormal bundle $\Cal I/\Cal I^2$ corresponds to
 $\Lambda ^4V$, i.e., to $\omega_4$ of SL(5); in
 particular, it is an irreducible $G$-bundle. But the
 nth symmetric power of $\Lambda ^4V$ is still an
 irreducible SL(5)-module, whence $S^n(\Cal I/\Cal
 I^2) \cong \Cal I^n/\Cal I^{n+1}$ is also an
 irreducible $G$-bundle. It is now a straightforward
 matter to prove the assertions, using the usual Bott
 vanishing theorems.

  \vskip3ex
   \head  3.  Descending to a subvariety\endhead

  \par (3.1) Suppose $X\subset\Bbb P^N$  is
 smooth and projectively normal with property
 ($\ast$): when does a subvariety $Y\subset X$ inherit
 this property?

   \proclaim{Proposition 3.2} Let $X\subset\Bbb P^N$ be
  smooth and projectively normal, $k\geq 2$ an integer
 for which
   \roster
   \item"(a)" $\Gamma (\Cal I_X(k)) \twoheadrightarrow
  \Gamma (\Cal
   I_X/\Cal I^2_X(k))$ is surjective.\endroster
 Let $Y\subset X$ be a subvariety, with ideal sheaf
  $\Cal J$.
 Assume that \roster
   \item"(b)" $H^1(X,\Cal J (1))=0$
   \item"(c)" $H^1(X,\Cal J^2(k))=0$
   \item"(d)" $H^1(X,\Cal I_X/\Cal I_X^2(k)\otimes \Cal
  J)=0$. \endroster
  Then for the embedding of $Y$ defined by the line
 bundle $\Cal O_Y(1)$,
   $$
   \Gamma (\Cal I_Y(k)) \twoheadrightarrow \Gamma
 (\Cal
  I_Y/\Cal I^2_Y(k)) \text{ is surjective.}
   $$ \endproclaim
   \flushpar {\it Proof.} Let $\Cal I'=\text{ Ker }(\Cal O_{\Bbb P^N}
 \to\Cal O_Y)$ be the ideal
 sheaf of $Y\subset\Bbb P^N$; then one has exact
 sequences
   $$\gather
   0\rightarrow\Cal I_X\rightarrow\Cal I'\rightarrow
  \Cal J\rightarrow 0\\
  \intertext{and (tensoring with $\Cal O_Y$ and noting
 $\Cal
  I_X/\Cal I^2_X$ is locally free)}
 0\rightarrow \Cal I_X/\Cal I_X^2\otimes\Cal O_Y
  \rightarrow \Cal I'/\Cal I^{\prime 2}\rightarrow \Cal J
  /\Cal J^2\rightarrow 0.
   \endgather$$
 Consider now for each $k$ the diagram
   $$\matrix\format \c&&\quad\c \\
   0 &\rightarrow &\Gamma(\Cal I_X(k)) &\rightarrow
  &\Gamma (\Cal I'(k))
  &\rightarrow &\Gamma (\Cal J(k)) &\rightarrow &0\\
  \vspace{2\jot}
 &&\sideset\and{\alpha_k}\to{\downarrow} \\
  \vspace{2\jot}
   && \Gamma (\Cal I_X/\Cal I^2_X(k))
  &&\sideset\and{\gamma_k}\to{\downarrow}
  &&\sideset\and{\delta_k}\to{\downarrow}  \\
  \vspace{2\jot}
 &&\sideset\and{\beta_k}\to{\downarrow} \\
  \vspace{2\jot}
   0 &\rightarrow &\Gamma(\Cal I_X/\Cal
  I^2_X(k)\otimes\Cal O_Y) &\rightarrow
  &\Gamma(\Cal I'/\Cal I^{\prime 2}(k)) &\rightarrow
  &\Gamma(\Cal J/\Cal J^2(k)) .
  \endmatrix
  $$
   The top row is exact by projective normality of
 $X\subset\Bbb P^N$. $\alpha _k$, $\beta_k$, and
  $\delta_k$ are
 surjective by hypotheses, hence so is $\gamma _k$.
 \par By (b), $\Cal O_Y (1)$ embeds $Y$ as a
 non-degenerate subvariety of its linear span
 $H\subset\Bbb P^N$;
   let $\Cal I_Y=\text{Ker
 }(\Cal
  O_H\rightarrow\Cal O_Y)$. By projection of $\Bbb
 P^N$ off a
  subspace onto $H$, one can split the surjection
 $\Cal I'\rightarrow\Cal I_Y\rightarrow 0$, deducing that
  $\Cal I_Y /\Cal I_Y^2$ is a direct summand of $\Cal
  I'/\Cal I^{\prime 2}$.  One has thus a surjection
  $$
   \Gamma(\Cal I'/\Cal I^{\prime 2}
  (k))\twoheadrightarrow\Gamma(\Cal I_Y/\Cal I_Y^2(k)).
  $$
 Composing with $\gamma_k $ gives a surjective map
  which factors through
 \flushpar $\Gamma(\Cal I_Y(k))
  \twoheadrightarrow\Gamma(\Cal I_Y/\Cal
  I^2_Y(k)),$ so this last map is also surjective.

   \proclaim{Corollary 3.3} Let $Y$ be a (possibly
  singular) projective variety, $A$ a very ample line
  bundle. Then there is an $m_0$ so that $m\geq m_0$
  implies $A^m$ satisfies condition {\rm ($\ast$)}; i.e., for
  the corresponding embeddings into projective space,
   $$
   H^1(\Cal I_Y^2 (k))=0,\qquad\text{ all } k\geq 3.
  $$\endproclaim
   \flushpar {\it Proof.} Fix the embedding
 $Y\subset\Bbb P^n=X$
   defined by $A$, with ideal sheaf $\Cal J$. Let
  $L=\Cal O_X(1)$, and $L'=L^m= \Cal O_X (m)$.  We
  claim that there is an integer $m_0$ so that for
  $m\geq m_0$, the embedding of $X$ into $\Bbb P^N$
  defined by $L'$ satisfies the conditions of Theorem
  3.2, for all $k\geq 3$. By Theorem 2.1, (a) is
  automatically true for all $m$ and all $k\geq 3$.
  (b) and (c) are certainly true for all $m\geq m_1$
  and all $k\geq 3$, since $\Cal J$ does not depend on
  $m$.

   To verify (d), denote by $\pi_i:X\times
 X\rightarrow X$ the projection maps $(i=1,2)$.
 Tensor the sequences
$$\gather
0 \rightarrow \Cal I_\Delta \rightarrow \Cal O_{X\times X} \rightarrow
\Cal O_\Delta \rightarrow 0\\
\intertext{and}
0 \rightarrow \Cal I_\Delta^2 \rightarrow
\Cal I_\Delta \rightarrow \Cal I_\Delta/\Cal I^2_\Delta \rightarrow 0
\endgather$$
with $\pi _2^\ast L'$, then apply $\pi _{1\ast}$; this yields
$$\gather
M_{L'}\cong \pi_{1\ast}(\pi^\ast_2 L'\otimes \Cal I_\Delta) \\
\Cal I_X/\Cal I_X^2(1)\cong\pi_{1\ast}(\pi^\ast_2 L'\otimes\Cal I_\Delta^2)
\endgather$$
for the embedding of $X$ given by $L'$. Thus
 $$\Cal I_X/\Cal I_X^2(k)\cong \pi_{1\ast}(\Cal
  I_\triangle^2\otimes(L^{'(k-1)}\boxtimes L')).
  $$
 Fix $E$, a locally free sheaf on $X$.  We claim
  there is an $m_2$ so that $m\geq m_2$ implies
 $$
  H^i(X,\Cal I_X/\Cal I_X^2(k)\otimes E)=0 \qquad
  \text{ for all }i>0,\ k\geq 2,\ L'=L^m.
  $$
 By the projection formula, it suffices to compute
 $$
  H^i (X,\pi_{1\ast}\{\Cal I_\triangle^2
 (L^{m(k-1)}\boxtimes L^m)\otimes \pi_1^\ast E\}).
  $$
  But the $R^q\pi_{1\ast}$ of the bracketed term
  vanishes for $q>0$ and $m$ large, so we may
 calculate
  the cohomology of the bracketed term on $X\times X$.
  Thus, it suffices to show that for $L$ ample and
  $\Cal F$ coherent on $X\times X$, there is a $j_0$ so that
  $$
  H^i (X\times X,\Cal F\otimes(L^j\boxtimes
  L^{j'})=0,\qquad
  \text{ all }i>0,j,j'\geq j_0.
  $$
  This is proved in the usual way by induction on $i$,
  writing $\Cal F$ as a quotient of a direct sum of
  terms of the form $L^s \boxtimes L^s$ and using the
  K\"unneth formulas plus standard vanishing theorems.

   To prove a comparable vanishing for $\Cal I_X/\Cal
  I_X^2(k)\otimes\Cal J$, take a projective resolution
  of the ideal sheaf $\Cal J$ by locally free $\Cal
  O_X$-modules.

 \proclaim{Corollary 3.4} Let $Y^n\subset\Bbb P^m$ be
 a (not necessarily smooth) complete intersection,
 defined by
  hypersurfaces of degrees $2\leq d_1\leq d_2\leq\dots
  \leq d_{m-n}$, where $n>0$. Consider the embedding
  of $Y$ defined by $\Cal O_Y (r)$, where
  $2r>d_{m-n}$. Then $H^1(\Cal I_Y^2  (k))=0, \ k\geq
  3$.\endproclaim
 \flushpar {\it Proof.} Embed
 $X=\Bbb
  P^m$ projectively normally via $L=\Cal O_{\Bbb P}(r)$ into some $\Bbb
P^N$.
  We verify the conditions of Proposition 3.2, noting the difference
between $L^k$ and $\Cal O_X (k)=\Cal O_{\Bbb P^m}(k)$. (a) is fulfilled
by Theorem 2.1, (b) is easy, and (c) is (1.2).
 For (d), we
  first deduce from (1.3.1) and from standard vanishing
  theorems on $X=\Bbb P^m$ that
  $$
   H^i(X,\Cal I_X/\Cal I_X^2(j))=0,\qquad 2\leq j<
  m,\text{ all }j;
  $$
 one must note that $H^2(X,\Cal I_X/\Cal I_X^2)=0$
  even though $H^1(X,\Omega^1_X)\neq 0$.  Thus, by
  the Koszul resolution of $\Cal J =\text{ Ker }(\Cal O_X\to \Cal O_Y)$,  there
is a surjection
  $$
   \oplus H^1(X,\Cal I_X/\Cal I_X^2\otimes
  L^k(-d_i))\twoheadrightarrow
  H^1(X,\Cal I_X/\Cal I_X^2\otimes L^k\otimes\Cal J).
  $$
 By (1.7.3), each summand on the left is of the form
 Coker $\Phi(L,L^{k-1}(-d_i))$.  Since the Gaussian of 2
 positive line bundles on $\Bbb P^m$ is surjective, we
 have the desired vanishing as long as $L^{k-1}
  (-d_i)$ is
 always positive, i.e., if $(k-1)r-d_i>0$.  In
  particular, this is true for all
   $k\geq 3$ and all $i$ once $2r>d_{m-n}$.
 \vskip1ex
 \par (3.5) One may argue more directly when $Y=X\cap H$ is a
 {\sl linear section} by a codimension $r$ subspace $H\subset \Bbb P^N$.
Thus, $Y$ is an l.c.i. (local complete intersection) in $X$ of
codimension $r$, and $\Gamma(\Cal O_X (1)\twoheadrightarrow
 \Gamma(\Cal O_Y(1)$.

 \proclaim{Proposition 3.6}  Let $X^n\subset\Bbb P^N$ be
 smooth, projectively normal, and arithmetically
 Cohen-Macaulay. Suppose $0<r<n$. Assume
 \vskip1ex
  \roster
  \item"(3.6.1)" The Gaussians $\Phi(\Cal O_X(1),\Cal O_X(k))$ are
  surjective, all $k\geq 1$. \endroster
 \vskip1ex
 \roster
  \item"(3.6.2)" $H^i (\Omega^1_X  (k-i))=0,\ 1\leq i
 \leq r-1,\ k\geq 2$.\endroster
 \vskip1ex
 \roster
  \item"(3.6.3)" $H^1(\Cal I_X^2(k))=0, \ k\geq 3.$
  \endroster
 \vskip1ex \flushpar
 Let $Y=X\cap H$ be a codimension $r$ linear section.
  Then $Y$ is projectively normal, and
 $$
  H^1(\Cal I_Y^2 (k))=0,\ k\geq 3.
 $$\endproclaim

 \flushpar {\it Proof}:  By hypothesis, $H^i (\Cal O
 _X(j)) =0,\ i\neq 0,n$.  A resolution of $\Cal J=$\newline
 Ker $(\Cal O_X\to\Cal O_Y)$ by $\Cal O_X$-modules is
 given by restriction of the Koszul complex defining $H$:
 $$
 0 \rightarrow \Cal O_X(-r)\rightarrow \dots \rightarrow
 \Cal O_X(-1)^r
 \rightarrow \Cal J \rightarrow 0.
 \tag3.6.4$$
 Thus, $H^1(\Cal J(k))=0$, all $k$,
  whence
 $$
 \Gamma (\Cal O_X(k)) \twoheadrightarrow \Gamma(\Cal O
 _Y(k)).
 $$
 It follows easily that $Y\subset H$ is projectively
 normal.
 \par Let $\Cal I_Y=\text{ Ker }(\Cal O_H
 \rightarrow\Cal O_Y).$ $Y$ is an l.c.i. in $X$,
 $H$, and $\Bbb P^N$, and $\Cal I_X/\Cal I_X^2$
 is locally free.  So,
 $$\gather
 \Cal I_X\otimes \Cal O_H \overset \sim\to\rightarrow
 \Cal I_Y \\
 \Cal I_X/\Cal I_X^2 \otimes \Cal O_Y
 \overset\sim\to{\rightarrow} \Cal I_Y /\Cal I_Y^2.
\endgather$$
 We thus have a commutative diagram with exact rows
   $$\gathered
  \CD  0 @>>> \Cal I_X \otimes \Cal J\\
  @.  @VVV\\
  0 @>>>  \Cal I_X/\Cal I^2_X \otimes \Cal J\endCD
  \CD {} @>>>  \Cal I_X \\
  @.  @VVV\\
  {} @>>>  \Cal I_X /\Cal I^2_X \endCD
  \CD  {} @>>>  \Cal I_Y\\
  @.  @VVV\\
  {} @>>>  \Cal I_Y /\Cal I^2_Y  \endCD
  \CD  {} @>>> 0\\
  @.  @. \\
  {} @>>> 0.\endCD
  \endgathered$$
 Taking cohomology and using (3.6.3), the conclusion would
 follow from

  $$\gather
  H^1(X, \Cal I_X/\Cal I_X^2 (k)\otimes \Cal J)=0, k\geq
 3.\\
\intertext{Via (3.6.4), it suffices to prove}
 H^i (X,\Cal I_X/\Cal I_X^2(k-i))=0,\ 1\leq i\leq r, \ k\geq 3.
 \tag3.6.5
 \endgather$$

 \par For $i=1$, this condition is
 $$
 H^1 (\Cal I_X/\Cal I_X^2 (k))=0,\ k\geq 2,
 $$
 which is a consequence of (3.6.1) via (1.7.3), (1.7.1),
 and the vanishing of $H^1(\Cal O_X(j))$.
 \par  For $i\geq 2$, compute using the conormal sequence
 of $X\subset \Bbb P^N$ and the Euler sequence for $\Omega
 ^1_{\Bbb P}$, restricted to $X$.  Then (3.6.5) is a
 consequence of (3.6.2).

 \proclaim{Corollary 3.7}  Let $Z\subset \Bbb P^3$ be a
 0-dimensional, length 5, arithmetically Gorenstein
 subscheme. Then
 $$
 H^1(\Cal I_Z^2(k))=0, \ k\geq 3.$$
 \endproclaim

 \flushpar {\it Proof.} By (\cite{S}, 4.2), $Z$ is
 Pfaffian, with ideal sheaf satisfying
 $$
  0 \to \Cal O_{\Bbb P} (-5) \to \Cal O_{\Bbb P} (-3)^{5}\to
 \Cal O_{\Bbb P} (-2)^{5} \to \Cal I_Z \to 0;
 $$
 so, $H^1 (\Cal I_Z(k))=0,\ k\geq 2.$ $Z$ is a codimension 6 linear
 section of the Grassmannian $X=G(2,5)\subset\Bbb P^9$.
 Imitating the proof of
 Proposition (3.6), one may verify the desired property for $Z$ by
checking (3.6.1), (3.6.2),
 and (3.6.3) for $r=n(=6)$, plus also
 $$
  H^n(\Cal O_X(k-n)) =0, \qquad k\geq2.
 \tag3.7.1$$
 But (3.6.1) is in \cite{W4}, (3.6.3) is (2.11); (3.6.2) and (3.7.1)
 use standard calculations (note $K_X\cong \Cal O_X(-5)$).

  \vskip3ex
 \head 4. Products and scrolls \endhead

  (4.1) If $L_i$ and $L_i'$ are line
 bundles on projective
 varieties $X_i,\ i=1,2$, one has bundles
 $N=L_1\boxtimes L_2$ and $N'=L_1' \boxtimes L_2' $ on
 $X_1\times X_2$.  A simple computation shows
 $$
 \Cal R_1(N,N')=\Cal R_1(L_1 ,L_1')\otimes\Gamma
 (L_2)\otimes\Gamma (L_2' )
 \ +\ \Gamma(L_1)\otimes\Gamma(L_1')\otimes\Cal R_1(L_2
 ,L_2'),
 \tag4.1.1$$
 where the two spaces on the right have intersection
  $$
  \Cal R_1(L_1 ,L_1')\otimes\Cal R_1(L_2 ,L_2').
  $$
 If $L_1=L_1'$ and $L_2 =L_2'$, so that $N=N'$, one
  has
 $$\multline
  \Cal R^s_1(N,N)=\Cal R^s_1 (L_1 ,L_1 )\otimes
 S^2\Gamma(L_2) + S^2\Gamma (L_1)\otimes\Cal
 R^s_1(L_2,L_2)\\ +\wedge^2
 \Gamma (L_1)\otimes\wedge^2\Gamma(L_2).
 \endmultline\tag4.1.2$$
 This says that the quadratic equations of the
  embedding defined by $N$ come from the quadratic
  equations of $X_1$ and $X_2$ and from the quadratic
  equations of the Segre embedding (the third term).
 \vskip1ex

 (4.2) Returning to the general case, the
 Gaussian map on $X_1\times X_2$ maps into the space
  $$
  \Gamma(\Omega^1_{X_1}\otimes L_1\otimes L_1')
  \otimes\Gamma(L_2 \otimes L_2')
  \ \oplus \
  \Gamma (L_1 \otimes L_1')\otimes \Gamma
  (\Omega^1_{X_2} \otimes L_2 \otimes L_2').
  $$
  Therefore, $\Cal R_2(N,N')$ is spanned by three
  subspaces:
 $$\multline
  \Cal R_2(N,N')=
 \Cal R_2 (L_1,L_1')\otimes\Gamma (L_2
 \otimes L_2')+\Gamma(L_1 \otimes L_1')\otimes\Cal
  R_2(L_2 ,L_2') \\
  +\Cal R_1(L_1 ,L_1')\otimes\Cal R_1(L_2 ,L_2')
  \endmultline\tag4.2.1$$
  (the last space goes to 0 on each component of the
 image of the Gaussian).

  \proclaim{Proposition 4.3}  Suppose $L_i$ is
 normally
  presented on $X_i\ (i=1,2)$, and that
  $$
   \Cal R^s_2  (L_i,L_i)\otimes\Gamma
  (L_i^k)\twoheadrightarrow  \Cal R_2(L_i,L_i^{k+1})
  \text{ is surjective,} \ k\geq 1.
  $$
 Then the same properties hold for $N=L_1 \boxtimes
  L_2$ on $X_1\times X_2$.
   \endproclaim
 \flushpar {\it Proof.} It is clear that $N$ is very ample
  and projectively normal; normal presentation is an easy
  exercise using (4.1). We check surjectivity of
  $$
  \Cal R^s_2 (N,N)\otimes \Gamma(N^k)
  \twoheadrightarrow
  \Cal R_2(N,N^{k+1}),\qquad \text{ all } k\geq 1.
  $$
 But $\Cal R^s_1=\Cal R^s_2$ is given by (4.1.2),
  while (4.2.1)
 gives $\Cal R_2$; it suffices to check surjectivity
  for each
 of the three types of terms which appear:
  $$\align
  \Cal R^s_1 (L_1 ,L_1 )\otimes S^2\Gamma (L_2 )
  \otimes
   \Gamma (L_1^k )
  \otimes
  \Gamma (L_2^k)&\twoheadrightarrow
  \Cal R_2(L_1 ,L_1^{k+1})
  \otimes \Gamma (L_2^{k+2})\\ \vspace{2\jot}
 S^2\Gamma (L_1 )\otimes \Cal R^s_1 (L_2 ,L_2 )
  \otimes \Gamma (L_1^k)\otimes \Gamma(L_2^k)
  &\twoheadrightarrow
  \Gamma (L_1^{k+1}) \otimes\Cal R_2(L_2 ,L_2^{k+1}
  )\\ \vspace{2\jot}
   \wedge^2 \Gamma (L_1)\otimes \wedge^2\Gamma (L_2 )
  \otimes \Gamma (L_1^k)
  \otimes \Gamma (L_2^k) &\twoheadrightarrow
  \Cal R_1(L_1,L_1^{k+1})\otimes\Cal R_1(L_2
  ,L_2^{k+1}).
  \endalign$$
 The surjectivity of the first two maps follows
  easily from the hypotheses. The last map is
  surjective by (1.11.b).

 \proclaim{Proposition 4.4}  Let $X=\Bbb P^n\times \Bbb
  P^m\subset\Bbb P^N$ be the Segre embedding, and
  $Y=X\cap H$ a linear section of dimension $>0$. Then
  $Y$ is normally presented, and
 $$
 H^1(\Cal I_Y^2(k))=0, \qquad \text{ all } k\geq 3 .
   $$
   \endproclaim
  \flushpar {\it Proof.} We verify the hypotheses of
   Proposition 3.6. $X\subset\Bbb P^N$ is well-known to be
  projectively normal, arithmetically Cohen-Macaulay,
  and normally presented.  The surjectivity of the
  Gaussians on $X$ follows from that on each factor
  ([W3], 4.12). The vanishing of (3.6.2) is an
  exercise using the K\"unneth formula; note e.g. that
 $$
 H^i (\Bbb P^n,\Omega^1_{\Bbb P}(a))\otimes
 H^j(\Bbb P^m,\Cal O(a))=0\ \text{for}\ 0<i+j<m,\
 \text{unless}\ i=1, j=0, a=0.
 $$
   Finally, the vanishing $H^1(\Cal I_X^2(k))=0$ for
  $k\geq 3$ is guaranteed by Proposition 4.3.
 \vskip1ex
 \par (4.5)
   Let $e_1\geq e_2\geq\dots\geq e_d\geq
  0$ be integers, $\Cal E =\oplus\Cal O_{\Bbb P}(e_i)$
   a vector bundle on $\Bbb P^1$, and $\pi:\Bbb P(\Cal E
  )\rightarrow\Bbb P^1$ the corresponding $\Bbb
  P^{d-1}$-bundle.  Write $f=\sum e_i\geq 2$. The {\sl
  normal scroll} $X(e_1,e_2,\dots,e_d)$ is the image
  of $\Bbb P(\Cal E)$ in $\Bbb P^r$ given by $\Cal O
  _{\Bbb P(\Cal E)}(1)\equiv H$; here $r=f+d-1$.
  Assume that $\Bbb P(\Cal E)\cong X\subset\Bbb P^r$,
  i.e., all $e_i>0$.  Pic$(X)$ is generated by
  $H$ and
  $R$, a fibre of $\pi$.
 \par It is well-known that $X\subset \Bbb P^r$ is defined determinantly as
  $$ rk \biggl( \aligned
  &x_{10} \dots \ \ \  x_{1e_1-1} \ \ \ x_{20}\dots\quad
  \dots x_{de_d-1}\\
  &x_{11} \dots \ \ \ x_{1e_1}\ \ \ \ \ \ x_{21}\dots\quad
  \dots  x_{de_d}
  \endaligned\biggr) \leq 1,
  $$
   where the $x_{ij}$ are the $\sum (e_i+1)=f+d$
  coordinates in $\Bbb P^r$ (we use that $e_i>0$).
 But the Segre embedding of $\Bbb P^{1}\times
  \Bbb P^{f-1}$ is similarly defined, as
   $$
  r k
  \biggl(\aligned
  & z_{11} \ \ \ z_{12} \dots \ \  z_{1f}\\
  & z_{21} \ \ \ z_{22} \dots \ \  z_{2f}
  \endaligned \biggr)\leq 1.
  $$
  Therefore, $X\subset\Bbb P^r$ may be viewed as a
  linear section
  of the Segre embedding. From the inclusion $X\subset\Bbb
  P^1\times\Bbb P^{f-1}$ followed by the projection
  maps, one has that $H$ is the restriction of $\Cal O
  (1)\boxtimes\Cal O(1)$ and $R$ is the restriction
  of $\pi_1^\ast\Cal O(1)$. Proposition 4.4 implies

  \proclaim{Proposition 4.6} Let $X=X(e_1,\dots,
  e_d)\subset\Bbb P^r$
 be a smooth rational normal scroll.  Then
  $X\subset\Bbb P^r$ is normally presented, and
   $$
  H^1(\Cal I_X^2 (k))=0,\qquad k\geq 3 .
  $$
   \endproclaim

  \vskip3ex
 \head 5. Pentagonal and tetragonal curves\endhead
 \par (5.1) Recall the results of F. Schreyer (\cite{S}, 6.7). Let
  $C\subset\Bbb P^{g-1}$ be a smooth canonical curve of genus
  $\geq 7$ with a base-point free $g^1_5$, i.e., a
  {\sl pentagonal curve}. The union of the linear spans of
  divisors in the $g^1_5$ form a four-dimensional
  scroll $X(e_1,e_2,e_3,e_4)\subset\Bbb P^{g-1}$.
  Write $f=\sum e_i=g-4$, and
 assume again that $e_4>0$, so that $X=\Bbb P(\Cal
  E)$.  Note $\Cal O_C(H)=K_C$, while $\Cal O_C(R)$ is the
  $g^1_5$. Let $\Cal J=$\newline $\text{Ker }(\Cal O_X
  \rightarrow\Cal O_C)$. Then there are integers
  $a_i,b_i$, at least $-1$, with $a_i+b_i=f-2,\ i=1,2,\dots,5$, and
  $\sum a_i=2g-12$, so that a resolution of $\Cal J$ as
  $\Cal O_X$-module is given by
   \vskip1ex \flushpar (5.1.1)
   $$
   0\rightarrow
  \Cal O_X(-5H+(f-2)R)\rightarrow
  \oplus O_X(-3H+b_iR)\rightarrow
  \oplus\Cal O_X(-2H+a_iR)\rightarrow
  \Cal J\rightarrow 0.
  $$
 In fact, $C$ is the Pfaffian locus of a skew-symmetric map
 $$
 \gather
 E^\ast\otimes E^\ast \rightarrow \Cal O_X (-5H + (f-2)R),
 \tag5.1.2\\
 E = \overset{5}\to{\underset{i=1}\to{\oplus}} \Cal O_X
 (3H-b_iR).
\endgather$$

 \par (5.2)  Straightforward calculations give
 $$\gather
 H^i (\Cal O (jH+ kR) =
 \left\{
 \matrix \format \l&&\quad\l \\
 0 &j=-1,-2,-3\text{ or }i=2\\
  H^i(\Bbb P^1, S^j\Cal E (k)) &j\geq 0\\
 H^{4-i} (\Cal O (-4-j)H + (f-2-k)R) &j\leq -4 .
  \endmatrix \right.
 \tag5.2.1\\ \vspace{1\jot}
   H^0(\Cal J(jH+kR))=0\text{ for }j\leq 1.\tag
  5.2.2\\  \vspace{1\jot}
   H^1(\Cal J(jH+kR))=0\text{ for }j\geq 2,k\geq
  0, \text{ or } j=1,\ k\leq 1.\tag 5.2.3
  \endgather$$
  \proclaim{Theorem 5.3}  Let $C\subset
  X(e_1,e_2,e_3,e_4)$ be a
 pentagonal curve as above, with all $e_i, a_i,
  b_i>0$.
 Then for the canonical embedding $C\subset\Bbb
  P^{g-1}$, one has
  $$
   H^1(\Cal I^2_C(k))=0,\qquad k\geq 3 .
  $$
  \endproclaim

  \flushpar {\it Proof.} We consider the embedding
 $C\subset
  X\subset\Bbb P^{g-1}$, and verify the
  hypotheses of Proposition 3.2. (a) is given by Proposition 4.6, and (b) is
contained in (5.2.3).

   To check condition (d),
  $$
  H^1(\Cal I_X/\Cal I_X^2 (k)\otimes\Cal J)=0,\qquad
  \text{for }k\geq 3,
  \tag5.3.1$$
  one needs via (5.1.1) the three vanishings
  $$\spreadmatrixlines{1\jot} \matrix \format\l&&\quad\l\\
   H^1(\Cal I_X/\Cal I_X^2(jH+a_iR))=0, &j\geq
  1\\
   H^2(\Cal I_X/\Cal I_X^2(jH+b_iR))=0, &j\geq 0\\
   H^3(\Cal I_X/\Cal I_X^2(jH+(f-2)R))=0, &j\geq -2.
  \endmatrix\tag 5.3.2$$
 Recall the standard short exact sequences on $X$:
   $$\matrix\format \c&&\quad\c \\
  (5.3.3) \qquad &0 &\rightarrow &\Cal I_X/\Cal
  I_X^2  &\rightarrow &\Omega^1_{\Bbb P}\vert_{X}
  &\rightarrow &\Omega^1_X &\rightarrow &0\\
  \vspace{3\jot}
  (5.3.4) \qquad &0  &\rightarrow &\Omega _{\Bbb
  P}^1 \vert_{X} &\rightarrow
  &\Gamma (H) \otimes \Cal O_X(-H) &\rightarrow &\Cal
  O_X &\rightarrow &0\\
  \vspace{3\jot}
  (5.3.5) \qquad &0  &\rightarrow  &\pi^\ast
  \Omega_{\Bbb P^1}^1 &\rightarrow &\Omega_X^1
  &\rightarrow  &\Omega_{X/\Bbb P^1}^1 &\rightarrow &0\\
  \vspace{3\jot}
  (5.3.6) \qquad &0  &\rightarrow &\Omega_{X/\Bbb
  P^1}^1  &\rightarrow  &\pi^\ast \Cal E(-H)
  &\rightarrow  &\Cal O_X &\rightarrow  &0.\\
   \endmatrix\ \ \ \ \ \ \ \ \ \
  $$

 These 4 sequences and the vanishing consequences of
  (5.2) are used to attack (5.3.2).  One also needs
  a few other easily verified facts:

  \roster\widestnumber\item{(5.6.10)} \item"(5.3.7)" The multiplication map
  $$
   \mu (H,jH+kR) : \Gamma (H)\otimes \Gamma
  (jH+kR)\rightarrow \Gamma ((j+1)H+kR)
 $$
  \item"" is surjective for $j\geq 0$, $k\geq 0.$
  (Proof: push to $\Bbb P^1$.)
  \item"(5.3.8)" The following map is surjective for
  $j\geq 0,\ k>0$:
  $$
   \Gamma (S^j \Cal E \otimes \Cal E (k)) \rightarrow
  \Gamma (S^{j+1} \Cal E (k)).
  $$
  \item"(5.3.9)" The Gaussian map
  $$
  \Phi(H,jH+kR): \Cal R (H,jH+kR)\rightarrow\Gamma
  (\Omega^1_X ((j+1)H+kR))
  $$
  \item"" is surjective for $j\geq 0,\ k>0$. (Proof: push to
  $\Bbb P^1$.)
  \item"(5.3.10)" One has natural isomorphisms on
 $\Bbb
  P^1$:\endroster
  $$
   \Cal O_{\Bbb P^1} \cong
  R^1 \pi_\ast \Omega^1_{\Bbb P} \vert_X \cong
  R^1\pi_\ast\Omega^1_X \cong
  R^1\pi_\ast \Omega^1_{X/\Bbb P^1} .
  $$
  (5.3.7)-(5.3.9) give surjectivity of some maps on
  global sections in twists of the sequences (5.3.4),
  (5.3.6), and (5.3.3).  (5.3.10) is used to show that
  in (5.3.3),
  $$
  H^1(\Omega^1_{\Bbb P}\vert_X(kR)) \cong
  H^1(\Omega^1_X (kR)),\qquad \text{ all } k>0.
  $$
 {}From these ingredients, the desired vanishing of
  (5.3.2) follows in a straightforward manner.

   It remains to verify condition (c) of Proposition 3.2,
  namely
  $$
   H^1(X,\Cal J^2(k))=0, \qquad \text{all } k\geq 3.
  $$
   First, a straightforward calculation using (5.1.1)
   and the conormal sequence for $C\subset X$ gives
  $$\gather
   \chi (\Cal J(3H))= \chi (\Cal J/\Cal J^2(3H))=10g-35,\tag5.3.11\\
 \intertext{whence}
 \chi (\Cal J^2(3H))=0.
  \tag5.3.12\endgather$$
  One has a surjection from (5.1.1)
 $$
 \oplus \Cal O_C((k-2)H+a_iR)\rightarrow \Cal
  J/\Cal J^2(kH) \rightarrow 0;
  $$
 since $a_i>0$ and $H$ restricted to $C$ is $K_C$,
  one has that
   $$
   H^1(\Cal J /\Cal J^2(kH))=0,\qquad k\geq 3 .
  $$
 The vanishing of $H^i(\Cal J (kH))=0$ for $i>0,
  \ k\geq 3$
 (use (5.1.1)), gives that
  $$
   H^i(\Cal J^2(kH))=0,\qquad i\geq 2, k\geq 3 .
  $$
 Combining with (5.3.12) yields
  $$
   h^0(\Cal J^2(3H))=h^1(\Cal J^2(3H)).\tag 5.3.13
  $$
   A fibre $R\simeq\Bbb P^3$ of $\pi$ is a smooth divisor
 not containing $C$, so
 $$
 \text{Tor}_i (\Cal O_C,\Cal O_R )=0,\ i>0,
 $$
 and (5.1.1) restricted to $R$ gives a projective
 resolution of $\bar{\Cal J}\subset\Cal O_{\Bbb P^3}$,
 the ideal sheaf of a 0-dimensional subscheme $Z$.
 Since $\Cal J/\Cal J^2$ is locally free, one shows
 similarly that $\Cal J^2$ restricted to a fibre is
 $\bar{\Cal J}^2$.  A simple count shows
 $$
 h^0 (\bar{\Cal J}^2(3))=
 h^1 (\bar{\Cal J}^2(3));
 $$
 but the second space is $0$, by (\cite{S}, 4.2), and
 Corollary 3.7. Since this is true for any fibre,
 $\Gamma(\Cal J^2 (3H))=0$. By (5.3.13), we conclude
  $$
  H^1(\Cal J^2(3H))=0.\tag 5.3.14
  $$

   Finally, we consider for $k\geq 4$ the commutative
  diagram
  $$\CD
  \oplus\Gamma(\Cal O_X((k-2)H+a_iR)) @>>> \Gamma(\Cal
  J(kH)) \\  @VVV  @VVV\\
  \oplus\Gamma(\Cal O_C((k-2)H+a_iR )) @>>>\Gamma
 (\Cal J/\Cal J^2(kH)).
  \endCD\tag 5.3.15
  $$
   The top and left vertical maps are surjective for
  $k\geq 4$ by easy calculation using (5.1.1) and
  (5.2). Tensoring (5.1.1) with $\Cal O_C$ gives an
  exact sequence
  $$
  \oplus\Cal O_C ((k-3)H+b_iR)\rightarrow
  \oplus\Cal O_C ((k-2)H+a_iR )\rightarrow
  \Cal J/\Cal J^2 (kH) \rightarrow 0.
  $$
 Since $H^1$ of the first term is 0 for degree
  reasons ($b_i>0$), one concludes surjectivity in the bottom
 row of (5.3.15).  Thus, the right vertical map is
  surjective, whence
  $$
   H^1(\Cal J^2(kH))=0,\qquad k\geq 4.\tag 5.3.16
  $$
 This completes the proof of the Theorem.

 \proclaim{Corollary 5.4} For a general pentagonal
  curve $C$ of genus $g\geq 8$, the canonical
 embedding satisfies
  $$
  H^1( \Cal I^2_C(k)) =0, \qquad k\geq 3 .
  $$
 \endproclaim

 \flushpar {\it Proof.} For $g=8$, the general curve lies
 on a scroll with all $e_i=1$, 4 $a_i$'s =1, and one=0
 (\cite{S}, 7.1). But it is easy to check that all the
 needed vanishing arguments of Theorem 5.3 still hold in
 that case (basically because $\Gamma(\Cal E(-2))=0$).

 Each $g\geq 9$ determines ordered sets of integers $\{
 e_i \},\ \{ a_i \},\ \text{and }\{ b_i \}$ which satisfy
 the appropriate equalities, are positive, and are {\sl
 balanced} (the numbers in a set differ by at most 1).  We
 claim there is a smooth connected pentagonal curve with
 these invariants; Theorem 5.3 then implies the Corollary.
 Let $X=\Bbb P(\Cal E )$. Following (5.1.2), let
 $$
 E^\ast =
 \underset {i=1} \to {\overset {5}\to\oplus}
 \Cal O_X (-3H+b_iR),\ L=\Cal O_X(-5H+(f-2)R);
 $$
 we show that a general skew-symmetric map
 $$
 \Psi : E^\ast \otimes E^\ast\rightarrow L
 $$
drops rank by 3 along a smooth connected curve.  We apply
 the following  useful observation of R. Lazarsfeld, but
 postpone the proof.
 \proclaim{Lemma 5.4.1} Let $\Bbb F\to X$ be a globally
 generated rank $\,r$
  vector bundle on a variety $X$, and $Z\subset \Bbb F$ a
 closed subvariety, with
 $$
 \dim (\text{{\rm Sing}}\, Z) <r.
 $$
 Then for a general section $s\in\Gamma (\Bbb F)$, the set
 $$
 \{ x\in X \,\vert\, s(x)\in Z\}
 $$
 is either empty, or smooth of dimension $\dim(Z)-r$.
 \endproclaim

 Returning to the proof of the Corollary, the map $\Psi$ is
 a section of
 $$
 F= \wedge^2 E \otimes L =
 \underset {j\neq k}\to{\oplus}
 \Cal O_X (H-(b_j-a_k)R).
 $$
 Since $\Cal O_X(H-kR)$ is globally generated if all
 $e_i\geq k$, one has global generation of $F$ if for all
 $i,j,k\ (j\neq k)$, one has
 $$
 e_i \geq b_j-a_k.
 \tag5.4.2$$
 It is an exercise to show that the balanced invariants of
 $g\geq 9$ satisfy (5.4.2), except when $g=10$ or 15
 (adding 5 to $g$ increases the maximum of $b_j-a_k$ by 1,
 while the minimum $e_i$ goes up by at least 1).  We
 exclude henceforth the special calculation needed to
 handle $g=10$ or 15.  So, apply Lemma 5.4.1 to the rank
 10 vector bundle associated to $F$, where $Z=$Pfaffian
 locus in $\Bbb F$ (where the rank of the skew-form drops
 by more than 1). Sing$(Z)$ is the 0-section of $\Bbb F$,
 hence has dimension 4.  We conclude that for a generic
 skew-form $\Psi$, either the rank drops by more than 1
 along a smooth curve, or nowhere.  The second case would
 give an exact sequence
 $$
 0\rightarrow A\rightarrow E^\ast \rightarrow E\otimes
 L\rightarrow B\rightarrow 0,
 $$
 where $A$ and $B$ are line bundles on $X$.  Restricting to
 a fibre of $\pi:\Bbb P(\Cal E)\to\Bbb P^1$ would
 give an exact sequence on $\Bbb P^3$ (cf. (5.1.1))
 $$
 0\rightarrow \Cal O (-5)\rightarrow \Cal O (-3)^{\oplus
 5}
 \rightarrow \Cal O (-2)^{\oplus 5} \rightarrow \Cal O
 \rightarrow 0,
 $$
 which is clearly impossible (apply $\Gamma$).  Thus, the
 cokernel of $E^\ast\rightarrow E\otimes L$ is the ideal
 sheaf $\Cal I_C$ of a smooth curve on $X$; its
 connectedness follows because $H^1(\Cal I_C)=0$,
 using (5.1.1).

  \vskip1ex
 \flushpar {\it Proof of Lemma 5.4.1}.  Let $\rho:X \times
 \Gamma (\Bbb F)\rightarrow \Bbb F$, a surjection of
 vector bundles; thus, $\bar{Z}\equiv \rho^{-1}(Z)$ and
 Sing$(\bar{Z})= \rho^{-1}(\text{Sing}(Z))$ have predicted
 dimensions.  Let $\pi :\bar{Z}\rightarrow \Gamma (\Bbb
 F)$ be the projection.  The hypothesis implies $\dim
 \text{ Sing }(\bar{Z})< \dim \Gamma (\Bbb F)$; so for
 some open $U\subset \Gamma (\Bbb F),\ \pi ^{-1}(U)\cap
 \text{ Sing }(\bar{Z})=\emptyset$.  Apply the theorem on
 generic smoothness to $\pi :\pi ^{-1}(U)\rightarrow U$.
 A fibre over $s\in \Gamma (\Bbb F)$ is
 $$
 \{ (x,s)\,\vert \, s(x)\in Z\} \cong s(X)\cap Z.
 $$
 If $\pi $ is dominant, the general fibre is smooth, of
 dimension
 $$
 \dim \bar{Z} - \dim \Gamma (\Bbb F) = \dim Z-r.
 $$
 If $\pi $ is not dominant, the general fibre is empty.
\vskip1ex

  (5.5) Better known
  than pentagonal curves are the tetragonal ones, with
  a base-point-free $g^1_4$ ([S], 6.2).  In this case,
  there is a 3-dimensional scroll $X$ defined by
  $e_i$'s which are positive unless the curve is
  bi-elliptic (i.e., a double cover of an elliptic curve);
  $C$ is the complete intersection on $X$ of divisors
  of the form $2H-b_iR$, where $0\leq b_1\leq b_2$,
  and $b_1+b_2=f-2=\sum e_i-2=g-5$.  Unless $C$ is
  bielliptic or lies on a del Pezzo surface, one has
  $b_1>0$.  An argument nearly identical to the one in
  Theorem 5.3 yields the following, whose proof is
  omitted:
   \proclaim{Theorem 5.6} Let $C$ be a tetragonal
  curve, with no $g^1_2$, $g^1_3$, or $g^2_5$, and which
  is not bielliptic (i.e., for which the scroll
  invariants $e_i>0$), with $g\geq 6$.  Then for the
  canonical embedding of $C$,
  $$\gather
   H^1(\Cal I^2_C (k))=0, \qquad k\geq 4\\
  \intertext{and}
  \dim H^1(\Cal I^2_C (3))=
  \left\{ \matrix\format\l&&\quad\l\\  g-7  &b_1>0\\  2(g-6)
  &b_1=0.\endmatrix \right.
  \endgather
  $$
  \endproclaim

  \proclaim{Corollary 5.7}  Let $C\subset\Bbb P^{g-1}$
  be a general canonical curve of genus $\geq 3$. Then
  $$
   H^1(\Cal I^2_C (k))=0, \qquad k\geq 3 .
  $$
 \endproclaim
 \flushpar {\it Proof.} Corollary 5.4 covers $g\geq
 8$;
  Theorem 5.6 covers $g=6$ and 7; and the general
  canonical curve of genus 3, 4, or 5 is a complete
  intersection, whence (1.2) applies.

 \vskip3ex
   \head 6. Vanishing for other canonical
  curves\endhead

  (6.1)  We are interested in
  determining which canonical curves $C\subset\Bbb
  P^{g-1}$ have the property
 $$
   H^1(\Cal I_C^2 (k))=0, \qquad k\geq 3.
 \tag"($\ast$)"
  $$
 It follows from Theorem 5.3 that the general
  pentagonal curve of genus $\geq 8$ has this
 property.
  On the other hand, by Theorem 5.6 this condition
  fails for every tetragonal curve with $g\geq 8$.  It
  seems reasonable to make the following

 \proclaim{Conjecture 6.2}
  Every canonical curve $C\subset\Bbb P^{g-1}$
  of Clifford index $\geq 3$ satisfies
  $$
   H^1(\Cal I^2_C (k))=0,\qquad k\geq 3 .
  $$
 \endproclaim

   \flushpar {\it Remarks} (6.3.1)  A hyperplane
  section of a $K3$ surface is a canonical curve, and
  any sufficiently positive section has the property
   ($\ast$), by (3.3).

  (6.3.2)  Any complete intersection curve in $\Bbb
  P^n$ satisfies $(\ast)$, except for the following
  cases (3.4): a plane curve of degree $\leq 6$; or a
  space curve with $d_1=2,\ d_2$=2, 3, or 4 (all of
  which have a $g^1_4$).

  (6.3.3)  According to a theorem of Schreyer and
  Voisin (e.g., [V]), a curve of Clifford index $\geq
  3$ is normally presented, with first syzygies
  generated by linear ones.  It follows easily (cf.
  [W2], 2.9.2) that $(T^2)_{-k}=0$ for $k\geq 4$.  By
  (1.6.2), one  concludes $H^1(\Cal I^2_C (k))=0$, all
  $k\geq 5$.  Conjecturally, Cliff(C)$\geq 4$ implies
  the syzygies are linear one further step; [W2] would
  then imply also that $H^1(\Cal I^2_C (4))=0$.

   (6.3.4)  This vanishing property is not as delicate
  as the question of the corank of the Gaussian
 $\Phi_K
  :\wedge^2\Gamma (K)\rightarrow \Gamma (K^{\otimes
  3})$, which is not surjective for any $K3$ curve [W2].

  \proclaim{Theorem 6.4}  Let $C\subset\Bbb P^{g-1}$
 be
  a general canonical curve, $3\leq g \leq 10,\ g\neq
  9$. Let $A$ be the affine cone over $C$.  Then
 $T^2_A=0$.\endproclaim

 \flushpar {\it Proof.}  For $3\leq g\leq 5$, $C$ (and
 hence $A$) is a complete intersection, so $T^2=0$; so
 assume $g\geq 6$.  By Corollary 1.6, one must show
 $H^1(\Cal
  I_C^2(k))=0$, all $k\geq 2$.  Now, vanishing for
  $k\geq 3$ follows from Corollary 5.10. By (1.8), it
  remains to show that the Gaussian $\Phi_K$ is
  injective for a general curve of genus 6, 7, 8, or
  10. But this result has been observed in [C-M].
 \vskip1ex

 \flushpar {\it Remarks} (6.5.1)  The above argument
  implies that for
 $C$ of genus 9 or $\geq 11$, $T^2$ must be non-zero
  in weight $-1$.
\par (6.5.2) Mukai proves \cite{Mk} that a general curve
 $C$ of genus 6, 7, 8, or 9 is cut out in a simple way
 from an appropriate homogeneous space $X_g=
 G/P\subset \Bbb P^{n(g)}$.  In the first 3 cases, one
 checks $\Cal I_X/\Cal I^2_X$ is irreducible
 (2.10), whence $H^1(\Cal I_X^2(k))=0$, all $k$.  As in
 \S3, one deduces $H^1(\Cal I_C^2(k))=0$, all $k$, $C$
 general, of genus 6, 7, or 8.  When $g=9$, $C$ is a
 linear section of $X_9 \subset \Bbb P^{13}$,
 corresponding to the weight $\omega_3$ of $Sp(3)=C_3$.
 The Gaussian of $X$ is surjective
  by [K], so a dimension count shows it has a
  one-dimensional kernel; since
 $$
   H^1(\Cal I_X^2 (2))\subset H^1(\Cal I_Y^2 (2))
  \ \text{for any linear section $Y$,}
  $$
  the Gaussian of a genus 9 curve cannot be injective
  (as proved in [C-M]).

 \vskip3ex

\head 7. Extendability of canonical curves \endhead
 \proclaim{Theorem 7.1}  Let $C\subset \Bbb P^{g-1}$
 be a canonical curve of genus $g\geq 8$, with
 {\rm Cliff}$(C)\geq 3$. Suppose
 $$
 H^1 (\Bbb P^{g-1},\Cal I^2_C(k))=0,\ k\geq 3.
 \tag"($\ast$)"
 $$
 Then $C$ is extendable if and only if the Gaussian
 $\Phi_K$ is not surjective.
 \endproclaim

 \flushpar {\it Proof.}  One implication is well-known
 \cite{W2}, even without ($\ast$); we prove the
 converse.  The affine ring of the cone is $A=\oplus
 \Gamma (C,K^{\otimes n}_C)= P/I$, where $P=\Bbb
 C[z_1,\dots,z_g]$; since Cliff$(C)\geq 3$, $I$ is
 generated by quadratic equations, with relations
 generated by linear ones (6.3.3).  By the
 non-surjectivity of $\Phi_K$, one has a non-trivial
 first-order deformation of weight -1 of $A$; we show
 it lifts to a deformation over $\Bbb C[t]$, i.e., a
 flat graded map
 $$
 \Bbb C[t] \rightarrow \Cal A = P[t]/\Cal I,
 $$
 where $t$ has degree 1.  Then $X=\text{Proj}\, (\Cal
 A)\subset \Bbb P^g$ is an extension of $X\cap\{
 t=0\}=C$; it is not the cone  over $C$ since the
 deformation is non-trivial to first-order.

 We follow the approach of \cite{St}, \S 2.  Write a
 minimal presentation of I:
 $$
 P^{\ell} \overset{r}\to\rightarrow P^k
 \overset{f}\to\rightarrow P \rightarrow
 P/I = A\rightarrow 0,
 $$
 where $r$ and $f$ correspond to matrices whose
 entries are homogeneous of degree respectively 1 and
 2.  The module of relation is
 $\Cal R=\text{Im}(r)\cong P^\ell/\text{Ker}(r)$; let
 $\Cal R_0\subset \Cal R$ be the submodule of
 \lq\lq trivial" (or Koszul) relations. One has
 $$
 T^2_A = \text{Coker} \{ r^t:\text{Hom}((P/I)^k,P/I)
 \rightarrow
 \text{Hom}(\Cal R/\Cal R_0,P/I)\},
 $$
 which we assume by ($\ast$) and (1.6) to be 0 in
 degree $\leq -2$. This means (cf \cite{W2}, (2.9.2)) that
 $$\align
 &\text{a map }\Cal R/\Cal R_0 \rightarrow A \text{
 sending homogeneous generators of }\Cal R \text{ to
 }\tag"($\ast\ast$)"\\
 &\text{homogeneous elements of degree }\leq 1\text{
 lifts to }(P/I)^k \rightarrow A.
 \endalign
  $$
 \par A deformation of weight $-1$ over $\Bbb
 C[\epsilon]/\epsilon^2$ is given by
 $$
 F=f+\epsilon f^{(1)},
 $$
 where $f^{(1)}:P^k\to P$ is a matrix of linear forms
 in $P$.  By flatness, there is a lifting of relations
 $$
 R=r+\epsilon r^{(1)},
 $$
 where $r^{(1)}:P^\ell \to P^k$, a matrix of
 constants, satisfies $F\cdot
 R=0\,\text{mod}\,\epsilon^2$. Since $f\cdot r=0$,
 $$
 f^{(1)}r + f r^{(1)} =0.
 \tag7.1.2$$
 This simply means that $f^{(1)}:P^k\to P$ vanishes
 mod $I$ on Im$(r)$, inducing a map $I\to P/I$.

 Try to lift the deformation and relation to $\Bbb
 C[\epsilon]/\epsilon^3$:
 $$\align
 F &= f+\epsilon f^{(1)} + \epsilon^2 f^{(2)}
 \tag7.1.3\\
 R &=r+\epsilon r^{(1)}.
 \endalign$$
 One seeks a matrix $f^{(2)}$ of constants; for
 homogeneity reasons, one could not perturb further
 the relations.  That $F\cdot
 R=0\,\text{mod}\,\epsilon^3$ requires
 $$
 f^{(1)}r^{(1)} + f^{(2)}r=0.
 \tag7.1.4$$
 Now, $f^{(1)}r^{(1)} : P^\ell\to P$ induces
 $$
 P^\ell \rightarrow P\rightarrow P/I,
 $$
 which we claim vanishes on Ker$( r )$: If
 $r(\alpha)=0$, then by (7.1.2) $f(r^{(1)}\alpha)=0$.
 As ker$(f)=\text{Im(r)}$, there's a $\beta \in
 P^\ell$ with
 $$
 r^{(1)}\alpha=r\beta .
 $$
 Thus,
 $$
 f^{(1)} r^{(1)} \alpha= f^{(1)}r\beta = -fr^{(1)}
 \beta \in I,
 $$
 verifying the claim.  Thus, $f^{(1)}r^{(1)}$ induces
 a homomorphism
 $$
 P^\ell /\text{Ker}(r) \cong \Cal R\rightarrow P/I.
 $$
 The map vanishes on $\Cal R_0$, as these relations
 involve the entries of $f$; one thus has a
 homomorphism
 $$
 \Cal R/\Cal R_0\rightarrow P/I,
 $$
 sending generators to elements of degree 1 (=deg of
 entries of $f^{(1)}r^{(1)}$).  By ($\ast\ast$), this
 lifts to a map
 $$
 g:(P/I)^k \rightarrow P/I,
 $$
 satisfying
 $$
 gr\equiv f^{(1)}r^{(1)} \text{ mod }I.
 $$
 Since these matrices have linear entries and $I$ is
 generated by quadrics, one may lift $-g$ to a map
 $$
 f^{(2)}:P^k\rightarrow P,
 $$
 satisfying (7.1.4).

 We assert that the equations in (7.1.3) satisfy
 $F\cdot R=0$ mod any power of $\epsilon$, i.e.,
 define a one-parameter flat deformation.  It suffices
 to show
 $$
 f^{(2)}r^{(1)} =0.
 $$

 As before, $f^{(2)}r^{(1)}$ induces a map
 $$
 P^\ell\rightarrow P/I,
 $$
  which we claim vanishes on Ker$(r)$:  If
 $r(\alpha)=0$, as above there is a $\beta \in P^\ell$
 with
 $$
 r^{(1)} (\alpha) =r(\beta).
 $$
 By (7.1.2) and (7.1.4), $0=f^{(2)}r(\alpha)=
 -f^{(1)}r^{(1)} (\alpha) = -f^{(1)}r(\beta)= fr^{(1)}(\beta)$,
 whence there is a $\gamma \in P^\ell$ with
 $$\gather
 r^{(1)} (\beta) =r(\gamma).\\
 \intertext{So,}
 f^{(2)}r^{(1)} (\alpha) = f^{(2)}r(\beta) =
 -f^{(1)}r^{(1)}(\beta) \text{ (by (7.1.4)) } =
 -f^{(1)}r(\gamma) = fr^{(1)}(\gamma)\in I.
 \endgather$$
 Again, $f^{(2)}r^{(1)}$ induces a map
 $$
 \Cal R/\Cal R_0\rightarrow P/I,
 $$
 sending generators to constants.  By hypothesis, this
 lifts to a map $g:P^k\to P/I$ satisfying
 $$
 gr\equiv f^{(2)} r^{(1)} \text{ mod }I.
 $$
 As the entries of $r$ have degree 1, while
 $f^{(2)}r^{(1)}$ has degree 0, one concludes that
 $$
 f^{(2)}r^{(1)} = 0,
 $$
 as desired to complete the proof.

 \vskip1ex

 \flushpar {\it Remark} 7.2:  One has actually proved
 a general statement about deformations:  Suppose
 $A=P/I$ is a graded $\Bbb C$-algebra, where $I$ is
 generated by quadrics, and the relations are
 generated by linear ones.  Suppose $T^2_A$ is 0 in
 degree $\leq -2$.  Then any first-order deformation
 of $A$ of negative weight lifts to a graded negative
 weight deformation over $\Bbb C[t]$.

 \vskip1ex

 \subhead Example 7.3 \endsubhead  Consider the
 canonical embedding of a smooth plane curve $C$ of
 degree $d\geq 7$.  Then ($\ast$) is satisfied,
 corank$( \Phi_K )=10$, $C$ sits on no K-3 surface,
 but $C$ is extendable to a singular rational surface
 with a non-smoothable simple elliptic singularity of
 degree $(3d-9)$.

 \flushpar {\it Proof.} ($\ast$) is given by Corollary 3.4,
 since $2(d-3)>d$. The corank of $\Phi$ is well known
 to be dim $\Gamma (-K_{\Bbb P})=10$.  That $C$ does
 not lie on a K-3 was proved in \cite{G-L}, but also
 follows from the deformation point of view \cite{W5}.
 By Theorem 7.1, $C\subset \Bbb P^{g-1}$ is
 extendable, which one sees explicitly as follows:
 \par Choose a smooth cubic $D$ intersecting $C$
 transversally in $Z$, a set of $3d$ points.  Consider
 the linear system $\vert \Cal I_Z(d)\vert$ of degree
 $d$ curves through $Z$; it contains $C$ and $D+E$,
 where $E$ is any degree $(d-3)$ curve, hence has
 dimension $g=$ \newline $(d-1)(d-2)/2$. One easily checks that
 the linear system has no base-points off $Z$, and off
 $D$ separates points and tangent vectors. It thus
 defines a morphism from $B=B\ell_Z(\Bbb P^2)$ into
 $\Bbb P^g$, with image a surface $X$ for which $C$ is
 a hyperplane section. Since $C$ is projectively
 normal, $X$ is normal and in fact is canonically
 trivial.  The proper transform $\widetilde{D}$ of $D$ in
 $B$ is an elliptic curve of self-intersection
 $-(3d-9)$, which gets collapsed to a point in $X$;
 this point is a simple elliptic singularity, known to
 be non-smoothable when the degree $3d-9$ is $>9$
 (\cite{P}).  In fact this non-smoothability for $d>6$
 is used in \cite{W5} to prove that $C$ is not in the
 closure of curves lying on a K-3 surface.

  \enddocument